\documentstyle[11pt]{article}
\newcommand{\R}{I\!\!R}
\newcommand{\N}{I\!\!N}
\newcommand{\C}{I\!\!\!\!C}

\newcommand{\Z}{{\bf Z}}
\newcommand\dint{\displaystyle\int}
\newcommand\dfrac{\displaystyle\frac}
\newcommand\dsum{\displaystyle\sum}
\newtheorem{thm}{Theora}[section]
\newtheorem{Th}[thm]{\footnotesize{\bf{Theorem}}}

\newtheorem{Cor}[thm]{\footnotesize{\bf{Corollary}}}
\newtheorem{Lem}[thm]{\footnotesize{\bf{Lemma}}}
\newtheorem{Prop}[thm]{\footnotesize{\bf{Proposition}}}
\newtheorem{Def}[thm]{\footnotesize{\bf{ Definition}}}
\newtheorem{Rem}[thm]{\footnotesize{\bf{Remark}}}
\newtheorem{Rems}[thm]{\footnotesize{\bf{Remarks}}}

\begin{document}
\title{\normalsize\bf{{Generalized Dunkl-Lipschitz Spaces}}}
\author{\small{\footnote{E.mail: Samir.Kallel@isimm.rnu.tn}\,
Samir Kallel
 }}
\date{}
\maketitle

\vspace*{-1cm}
\begin{center}
\footnote{Fax Number: 0021671885350.}{\footnotesize{Department of
Mathematics, Faculty of Sciences of Tunis,\\ University Campus,
2092 Tunis, Tunisia}}
\end{center}
\vspace*{1cm}
\begin{abstract}
This paper deals with generalized Lipschitz spaces $\wedge^k_{\alpha,p,q}(\R)$ in the context of Dunkl harmonic analysis on $\R$ , for all real $\alpha$. It also introduces a generalized Dunkl-Lipschitz spaces ${\cal T}\wedge^k_{\alpha,p,q}(\R^2_+)$ of
$k$-temperature on $\R^2_+$.
 Some properties and continuous embedding of these spaces and the isomorphism of ${\cal T}\wedge^k_{\alpha,p,q}(\R^2_+)$ and $\wedge^k_{\alpha,p,q}(\R)$ are established.
\end{abstract}
{\scriptsize{ {\bf{Keywords :}} Dunkl operator.  Poisson transform. Heat transform.
 Dunkl transform.  Generalized Dunkl-Lipschitz space.}}\vspace*{-0.3cm}
\section{\footnotesize{{\bf{Introduction}}}}
\hspace*{5mm} In \cite{S-K}, we have introduced and characterized for $\alpha>0$ and $1\leq p,q\leq\infty$ the generalized Dunkl-Lipschitz spaces $\wedge^k_{\alpha,p,q}
(\R)$ associated with the Dunkl operator with parameter $k\geq0$
$${\cal D}_kf(x)=f'(x)+k\dfrac{f(x)-f(-x)}{x},\,\,\,f\in C^1(\R).$$ We were interested in characterizing the functions $f\in\wedge^k_{\alpha,p,q}(\R)$ for $\alpha>0$ in terms of their $k$-Poisson transform and the second order $L^p_k$-modulus of continuity. It is natural to extend the theory of the spaces $\wedge^k_{\alpha,p,q}(\R)$ for all real $\alpha$. To get this extension we use the $k$-heat transforms, since it is better suited in the treatment of tempered distributions than the $k$-Poisson transforms. More precisely, we define the spaces $\wedge^k_{\alpha,p,q}(\R)$ for $\alpha\leq0$ as spaces of tempered
distributions $T$ that belongs to an appropriate Lebesgue space for which the $k$-heat
transform $G^k_t(T)$ of $T$ satisfies the condition
 $$
\left\{\dint_0^1t^{q(n-\frac{1}{2}\alpha)}\|\partial^n_tG^k_t(T)\|_{k,p}^q\dfrac{dt}{t}\right\}^{\frac{1}{q}}<\infty,\;\;\;\;\;\;\;\mbox{if}\;\;\;1\leq
q<\infty$$ and $$\;\;\;\;\;\;\;\;\; \displaystyle\sup_{0<t\leq1}t^{n-\frac{1}{2}\alpha}\|\partial^n_tG^k_t(T)\|_{k,p}<\infty,\;\;\;\;\;\;\mbox{if}\;\;\;q=\infty,$$
where $n=\overline{(\frac{\alpha}{2})}$ and $\overline{\alpha}$ is the smallest non-negative integer larger than $\alpha$. The first goal of this paper is to study these spaces. As it is well known, the fractional integral operators play an important role in this theory. Here we use the Dunkl-Bessel potential ${\cal J}^k_{\alpha}$ which we show that ${\cal J}^k_{\alpha}$ is a topological isomorphism from $\wedge^k_{\alpha,p,q}(\R)$ onto $\wedge^k_{\alpha+\beta,p,q}(\R)$, with $1\leq p,q\leq\infty$ and $\alpha$, $\beta\in\R$. Next, certain properties and continuous embedding for $\wedge^k_{\alpha,p,q}(\R)$ are given.\\ Our second objective will study the generalized Dunkl-Lipschitz spaces of $k$-temperatures (i.e., solutions of the Dunkl-type heat equation
 $({\cal D}^2_k-\partial_t){\cal U}=0$) on the whole half-plane $\R^2_+=\left\{(x,t):x\in\R,t>0\right\}$ which denote by ${\cal T}\wedge^k_{\alpha,p,q}(\R^2_+)$, $1\leq p,q\leq\infty$. In Theorem \ref{Th.VI'.1}, we prove some basic properties of the space ${\cal T}\wedge^k_{\alpha,p,q}(\R^2_+)$ in which the most important one is the fact that the topological property of the space ${\cal T}\wedge^k_{\alpha,p,q}(\R^2_+)$ does not depend on the (Lipschitz) index $\alpha$. Thus, we should ask what relations there are between the generalized Dunkl-Lipschitz spaces $\wedge^k_{\alpha,p,q}(\R)$ and the generalized Dunkl-Lipschitz spaces of $k$-temperatures ${\cal T}\wedge^k_{\alpha,p,q}(\R^2_+)$. To reply to this question we must use the $k$-heat transforms. In Theorem \ref{Th.VII.1}, we establish that a $k$-temperature ${\cal U}$ belongs to ${\cal T}\wedge^k_{\alpha,p,q}(\R^2_+)$ if and only if it is the $k$-heat transform of an element of $\wedge^k_{\alpha,p,q}(\R)$. So that the spaces $\wedge^k_{\alpha,p,q}(\R)$ for $\alpha\leq0$, which consist of tempered distributions, can be realized as spaces of functions.\par Similar results have been obtained by T. M. Fleet and M. H. Taibleson \cite{Flett1,Taibleson} in the framework of classical case $k=0$. Later, R. Johnson \cite{Johnson}, adopting Flett's idea, defined a space of temperatures which is isomorphic to the Lipschitz space of Herz. His method leaned on a theory of Riesz potentials for temperatures. Additionally, for $\alpha>0$, the generalized Dunkl-Lipschitz spaces or Besov-Dunkl spaces have been studied extensively by several mathematicians and characterized in different ways by many authors (see \cite{Ch-An-Sa-Si, Ch-Sa, Ch-Si, BLA, S-K, Lot1,Lot2}).
 \par In this work, it is important to mention that the $1$D restriction is due to the fact that Dunkl translations operations in higher dimension are not yet known to be bounded on $L^p_k$ apart from $p=2$.
 \par The organization of this paper is as follows. In Section 2, we recall some basic harmonic analysis results related to Dunkl operator. In Section 3, we recall some properties of the $k$-heat transform of a measurable function. In Section 4, a semi-group formula for $k$-temperatures is proved which will be used frequently. In Section 5, the Dunkl-Bessel potential is defined and related properties are investigated. In Section 6, $\wedge^k_{\alpha,p,q}(\R)$ for real $\alpha$ is defined and its properties have been obtained. In this section we also proved that ${\cal J}^k_{\beta}$ is a topological isomorphism from $\wedge^k_{\alpha,p,q}(\R)$ onto $\wedge^k_{\alpha+\beta,p,q}(\R)$, $\alpha,\beta\in\R$, and a variety of equivalent norms for $\wedge^k_{\alpha,p,q}(\R)$ are given. The remainder of this section is devoted to some properties and continuous embedding for $\wedge^k_{\alpha,p,q}(\R)$. In section 7, we defined the space ${\cal T}\wedge^k_{\alpha,p,q}(\R^2_+)$, the equivalence of several norms on ${\cal T}\wedge^k_{\alpha,p,q}(\R^2_+)$ is proved and some properties of this space are studied. At the end, the isomorphism of ${\cal T}\wedge^k_{\alpha,p,q}(\R^2_+)$ and $\wedge^k_{\alpha,p,q}(\R)$ is established.
\par In what follows, $B$ represents a suitable positive
constant which is not necessarily the same in each occurrence.
\section{\footnotesize{{\bf{Preliminaries in the Dunkl Setting on $\R$}}}}\par
\hspace*{5mm} In this section we state some definitions and
results which are useful in the sequel and we refer for more details to the articles
 \cite{N-S,Dunkl3,Dunkl1,Dunkl2}, \cite{Jeu1}, \cite{Xu1} and \cite{Trim2}. We first begin by some notations.\\[4mm]
{\bf{Notations}} \begin{itemize} \item $C_0(\R)$ is the space of
continuous functions vanishing at infinity, equipped with the
usual topology of uniform convergence on $\R$. \item ${\cal
E}(\R)$ is the space of $C^{\infty}$-functions on $\R$, endowed with
the usual topology of uniform convergence of any derivative on compact subsets of $\R$.
\item $S(\R)$ is the space of $C^{\infty}$-functions on $\R$ which
are rapidly decreasing as well as their derivatives, endowed with the
topology defined by the semi-norms
$$\rho_{s,l}(\varphi):=\displaystyle\sup_{x\in\R,j\leq s}(1+x^2)^l|{\cal D}_k^j
\varphi(x)|,\,\,s,l\in\N.$$  \item $S'(\R)$ is the space of tempered
distributions on $\R$ which is the topological dual of $S(\R)$.
\end{itemize}
\par  The Dunkl operator ${\cal D}_k$ with parameter $k\geq0$ is given by
$${\cal D}_kf(x) :=f'(x)+k\dfrac{f(x)-f(-x)}{x}\, ,\quad f\in
C^1(\R).$$ For $k=0$, ${\cal D}_0$ reduces to the usual derivative
which will be denoted by ${\cal D}$. The Dunkl intertwining
operator $V_k$ is defined in \cite{Dunkl1} on polynomials $f$ by
$$
{\cal D}_kV_kf=V_k{\cal D}f\,\,\mbox{and}\,\,V_k1=1.
$$
For $k>0,\,V_k$ has the following representation (see
\cite{Dunkl1}, Theorem 5.1)
\begin{equation}\label{e:I.1}
V_kf(x):=\dfrac{2^{-2k}\Gamma(2k+1)}{\Gamma(k)\Gamma(k+1)}\dint_{-1}^1
f(xt)(1-t^2)^{k-1}(1+t)dt.
\end{equation}
This integral transform extends to a topological automorphism to the space
${\cal E}(\R)$ (see \cite{Trim2} and \cite{N-S}). For $k\geq 0$, and
$\lambda\in\C$, the initial problem $$\left\{
\begin{array}{rcl}
{\cal D}_ku(x)&=&\lambda u(x),\,x\in\R,\\ u(0)&=&1,
\end{array}
\right.$$ has a unique analytic solution $u(x)=E_k(\lambda, x)$, called
Dunkl kernel \cite{Dunkl1} and given by
$$E_k(\lambda, x):=j_{k-\frac{1}{2}}(i\lambda
x)+\dfrac{\lambda x}{2k+1}j_{k+\frac{1}{2}}(i\lambda x),$$ where
$j_{\alpha}$ is the normalized Bessel function, defined
for $\alpha\geq-\dfrac{1}{2}$ by
$$j_{\alpha}(z):=\Gamma(\alpha+1)\dsum_{n=0}^{+\infty}\dfrac{(-1)^n}{n!}
\dfrac{(\frac{z}{2})^{2n}}{\Gamma(n+\alpha+1)} ,\,z\in\C.$$ We
remark that $E_k(\lambda, x)=V_k(e^{\lambda .})(x)$. Formula
(\ref{e:I.1}) and the last result imply that
\begin{equation}\label{e:I.2}
\mid E_k(\lambda, x)\mid\leq e^{\mid\lambda\mid\mid x\mid},\,\mid
E_k(\lambda,x)\mid\leq e^{\mid x\mid\mid{\cal R
}e\lambda\mid},\,\mid E_k(-iy,x)\mid\leq 1,
\end{equation}
for all $x,y\in\R$ and $\lambda \in\C$.
\\ For all $f$ and $g$ in $C^1(\R)$ with at least one of them is even, we have
$$
{\cal D}_k(fg)=({\cal D}_kf)g+g({\cal D}_kg).
$$
For $f\in C^1_b(\R)$ and $g$ in $S(\R)$, we have
$$
\dint_{\R}{\cal D}_kf(x)g(x)|x|^{2k}dx=-\dint_{\R}f(x){\cal
D}_kg(x)|x|^{2k}dx.
$$
 Hereafter, we denote by $L^p(\R,|x|^{2k}dx)$, $p\in[1,\infty]$, the space of measurable functions on $\R$ such that
$$\|f\|_{k,p}:=(\dint_{\R}|f(x)|^p|x|^{2k}dx)^{\frac{1}{p}}<+\infty,\,\,\,\mbox\,\,\,
1\leq p<\infty,$$ and
$$\|f\|_{k,\infty}:=ess\!\displaystyle\sup_{x\in\R}|f(x)|<+\infty.$$
\par The Dunkl kernel gives rise to an integral transform, called
Dunkl transform on $\R$, which was introduced by Dunkl in
\cite{Dunkl2}, where already many basic properties were
established. Dunkl's results were completed and extended later on
by de Jeu in \cite{Jeu1}. \par The Dunkl transform of a function
$f\in L^1(\R, |x|^{2k}dx)$ is given by
$$\forall y\in\R,\,\,\,{\cal F}_k(f)(y):=c_k\dint_{\R}f(x)E_k(x,-iy)|x|^{2k}dx,$$
where $c_k:=\frac{1}{2^{k+\frac{1}{2}}\Gamma(k+\frac{1}{2})}$.
\par We summarize the properties of ${\cal F}_k(f)$ in the following proposition :
\begin{Prop} \cite{Jeu1}\label{Prop.I.1}
\par (i) For all $f\in S(\R)$, we have
$${\cal F}_k({\cal D}_kf)(x)=ix{\cal F}_k(f)(x),\,\,\,x\in\R.$$
\quad (ii) {\footnotesize\bf{Inversion formula :}} For all
$f\in L^1(\R,|x|^{2k}dx)$ such that  ${\cal F}_k(f)$ belongs to
$L^1(\R,|x|^{2k}dx)$, we have
$$f(x)=\dint_{\R}E_k(x,iy){\cal F}_k(f)(y)|y|^{2k}dy\,\,\,\,a.e.$$
\quad (iii) {\footnotesize\bf{Plancherel's Theorem :}} The
Dunkl transform extends to an isometry of $L^2(\R,|x|^{2k}dx)$. In
particular, we have the following Plancherel's formula
$$\|f\|_{k,2}=\|{\cal F}_k(f)\|_{k,2},\,\,\,f\in L^2(\R,|x|^{2k}dx).$$
\end{Prop}
\begin{Def}
Let $f\in C(\R)$ (denotes the space of continuous functions on $\R$) and
$y\in\R$. Then ${\cal T}^k_yf(x)=u(x,y)$ is
the unique solution of the following Cauchy problem
$$\left\{
\begin{array}{rcl}
{\cal D}_{k,x}u(x,y)&=&{\cal D}_{k,y}u(x,y),\\ u(x,0)&=&f(x).
\end{array}
\right.$$ ${\cal T}^k_y$ is called the Dunkl translation
operator.
\end{Def}
\begin{Rem}\label{Rem.I.1} In what follows we point out some remarks.\\
\begin{itemize}
 \item The operator ${\cal T}^k_x$ admits the following integral
 representation
\begin{equation}\label{e.I.3}
{\cal
T}^k_yf(x):=d_k\left(\dint_0^{\pi}f_e(G(x,y,\theta))h^e(x,y,\theta)\sin^{2k-1}\theta
d \theta \right.\hspace*{6cm}
\end{equation}
$$\hspace*{3cm} \left. +\dint_0^{\pi}f_o(G(x,y,\theta))h^o(x,y,\theta)\sin^{2k-1}\theta
d\theta \right),$$ where
$$d_k:=\dfrac{\Gamma(k+\frac{1}{2})}{\Gamma(k)\Gamma(\frac{1}{2})},\,\,\,G(x,y,\theta)=\sqrt{x^2+y^2-2|xy|\cos\theta},\,\,\,h^e(x,y,\theta)=
1-sgn(xy)\cos\theta,$$
$$h^o(x,y,\theta)=\left\{
\begin{array}{rcl}
\dfrac{(x+y)h^e(x,y,\theta)}{G(x,y,\theta)}&,&\mbox{if}\,\,\,(x,y)\neq(0,0),\\
0\hspace*{1cm}&,&\mbox{otherwise},
\end{array}
\right.$$
$$f_e(x)=\dfrac{1}{2}(f(x)+f(-x))\,\,\,\mbox{and}\,\,\,f_o(x)=\dfrac{1}{2}(f(x)-f(-x)).$$
\item There is an abstract formula for ${\cal T}^k_y$, $y\in\R$, given in terms of the intertwining operator $V_k$ and its inverse, ( see \cite{Trim2,N-S} ).
It takes the form of
$${\cal T}^k_yf(x):=(V_k)_x\otimes(V_k)_y\left[(V_k)^{-1}(f)(x+y)\right],\,\,\,\,x\in\R,\,\,f\in{\cal E}(\R).$$
\item The Dunkl translation operators satisfy for $x,y\in\R$ the
following relations
$$
\begin{array}{rcl}
{\cal T}_x^kf(y)={\cal T}_y^kf(x)\quad ,\quad {\cal T}_0^kf(y)=f(y),\\
{\cal T}_x^k{\cal T}_y^k={\cal T}_y^k{\cal T}_x^k\quad ,\quad
{\cal T}_x^k{\cal D}_k={\cal D}_k{\cal T}_x^k.
\end{array}$$
\item For each $y\in\R$, the Dunkl translation operator ${\cal T}^k_y$ extends to a bounded operator on
$L^p(\R,|x|^{2k}dx)$. More precisely
\begin{equation}\label{e.I.4}
\|{\cal T}^k_yf\|_{k,p}\leq 3\|f\|_{k,p},\,\, 1\leq p\leq\infty.
\end{equation}
\item Unusually, ${\cal T}^k_y$ is not a positive operator in general (see \cite{Ros}), but if $f$ is even, then ${\cal T}^k_yf(x)=d_k\int_0^\pi f(G(x,y,\theta))h^e(x,y,\theta)\sin^{2k-1}\theta d\theta$, which shows that ${\cal T}^k_yf(x)\geq0$ whenever $f$ is non-negative.
\item From the generalized Taylor formula with integral remainder (see \cite{Mourou}, Theorem 2 p. 349), we have for $f\in {\cal E}(\R)$ and $x,y\in\R$
\begin{equation}\label{e.II.Taylor}
\left({\cal T}^k_xf-f\right)(y)=\dint_{-|x|}^{|x|}\left(\dfrac{sgn(x)}{2|x|^{2k}}-\dfrac{sgn(z)}{2|z|^{2k}}\right){\cal T}^k_z({\cal D}_kf)(y)|z|^{2k}dz.
\end{equation}
\end{itemize}
\end{Rem}
 \par Associated to the Dunkl translation operator ${\cal T}^k_y$, the Dunkl convolution $f\ast_k g$ of two appropriate functions $f$ and $g$ on
$\R$ defined by
$$f\ast_k g(x):=\dint_{\R}{\cal T}^k_xf(-y)g(y)|y|^{2k}dy,\,\,\,x\in\R.$$ The Dunkl convolution preserves the main properties of the classical convolution which corresponds to $k=0$.\par For $S\in S'(\R)$ and $f\in S(\R)$, we
define the Dunkl convolution product $S\ast_kf$ by
$$S\ast_kf(x):=<S_y,{\cal T}^k_xf(-y)>.$$

\section{\footnotesize{{\bf{The $k$-Heat
Transforms of a Function}}}} \par We recall
some properties of the $k$-heat transforms of a measurable
function $f$ and we refer for more details to the survey \cite{N-A-S} and the references therein.
\par - For $t>0$, let $F^k_t$ be the function defined by
$$F^k_t(x):=(2t)^{-(k+\frac{1}{2})}e^{-\frac{x^2}{4t}}$$
which is a solution of the Dunkl-type heat equation $({\cal D}^2_k-\partial_t){\cal U}=0$ on the half-plane $\R^2_+$, \footnote{$\R^2_+=\left\{(x,t):x\in\R,t>0\right\}$}. The function $F^k_t$ may be
called the heat kernel associated with Dunkl operator or the $k$-heat kernel and it has the following basic properties :
\begin{Prop}\label{Prop.II.1}
For all $t>0$ and $n,\,m\,\in\N$, we have
\\ (i) ${\cal F}_k(F^k_t)(x)=e^{-tx^2}$ and
$\int_{\R}F^k_t(x)|x|^{2k}dx=c_k^{-1}$.
\\ (ii) $\int_{\R}|{\cal D}_k^nF^k_t(x)||x|^{2k}dx\leq
B(k,n)t^{-\frac{n}{2}}$.
\\ (iii)
$\partial^m_tF^k_t(x)=t^{-m}R(\frac{x^2}{4t})F^k_t(x)$, where $R$
is a polynomial of degree $m$ with coefficients depending only on
$m$ and $k$.
\\ (iv) $\int_{\R}|\partial^m_tF^k_t(x)||x|^{2k}dx\leq
B(k,m)t^{-m}$ and $\int_{\R}\partial^m_tF^k_t(x)|x|^{2k}dx=0$.
\end{Prop}
\begin{Def}
The $k$-heat transform of a smooth measurable function $f$ on $\R$ is given by
$$G^k_t(f)(x):=F^k_t\ast_kf(x),\,\,t>0.$$
\end{Def}
\begin{Th}\label{Th.II.1}
Let $f$ be a measurable bounded function on $\R$. Then,
\\ (i) $(x,t)\longmapsto G^k_t(f)(x)$ is infinitely differentiable on
$\R^2_+$ and it is a
solution of the Dunkl-type heat equation. Further, if
$n,m\in\N$, then for all $t>0$
$${\cal D}^n_kG^k_t(f)={\cal D}^n_kF^k_t\ast_kf\,\,\,\mbox{and}\,\,\,\partial^n_tG^k_t(f)=\partial^n_tF^k_t\ast_kf.$$
(ii) For all $s,t>0$ and $x\in\R$, we have
$G^k_{t+s}(f)(x)=\int_{\R}{\cal
T}^k_{-y}F^k_t(x)G^k_s(f)(y)|y|^{2k}dy$.
\\ (iii) If $f\in C_b(\R)$, then
$G^k_t(f)(x)\longrightarrow f(\xi)$ as
$(x,t)\longrightarrow(\xi,0)$.
\end{Th}
\begin{Th}\label{Th.II.2}
Let $p\in[1,\infty]$ and let $f\in L^p(\R,|x|^{2k}dx)$. Then the
$k$-heat transform $G^k_t(f)$ of $f$ has the following properties
:
\\ (i) For all $t>0$ and $m\in\N$, we have
 $$\|G^k_t(f)\|_{k,p}\leq c_k^{-1}\|f\|_{k,p}\,\,\,\mbox{and}\,\,\,\|\partial^m_tG^k_t(f)\|_{k,p}\leq B(k,m)t^{-m}\|f\|_{k,p}.$$
 (ii) If $1\leq p<r<\infty$ and
 $\delta=\frac{1}{p}-\frac{1}{r}$, then for all $t>0$
 $$\|G^k_t(f)\|_{k,r}\leq
 t^{-(k+\frac{1}{2})\delta}c_k^{\delta-2}\|f\|_{k,p}$$ and $\|G^k_t(f)\|_{k,r}=\circ(t^{-(k+\frac{1}{2})\delta})$, \footnote{
    $f(x)=\circ(g(x))$, $x\longrightarrow a$, means $f(x)/g(x)\longrightarrow0$ as $x\longrightarrow a$. }, as $t\longrightarrow0^+$.
\end{Th}
\begin{Def}
For any $T\in S'(\R)$, the $k$-heat transform of $T$ is given by
$$G^k_t(T)(x):=T\ast_kF^k_t(x),\,\,x\in\R.$$
\end{Def}
\section{\footnotesize{{\bf{A Semi-group Formula for
$k$-Temperatures}}}} \hspace*{5mm} Hereafter we shall be concerned
 mostly with temperatures associated with the Dunkl setting on $\R$ which we recall the $k$-temperatures, satisfying a property which we call
 "semi-group formula".
 \begin{Def}
A function ${\cal U}$ on $\R^2_+$ is said to be a $k$-temperature if it
is indefinitely differentiable on $\R^2_+$ and satisfies at each
point of $\R^2_+$ the Dunkl-type heat equation i.e.,
$${\cal D}^2_k{\cal U}(x,t)=\partial_t{\cal U}(x,t).$$
\end{Def}
\par - We consider the following initial-value problem for the $k$-heat equation :
$$(IVP)\left\{
\begin{array}{rcl}
({\cal D}_k^2-\partial_t){\cal U}=&0&\mbox{on}\,\,\,\R^2_+\\
{\cal U}(.,0)=&f&
\end{array}
\right.$$ with initial data $f\in C_b(\R)$ ( that is, the space of bounded continuous functions on $\R$). For $f\in C_0(\R)$, the function
$$H_tf(x)=\dint_{\R}{\cal T}^k_{-y}F^k_t(x)f(y)|y|^{2k}dy,\,\,t>0,$$
solves initial value problem (IVP) (see \cite{Ros-Voit}).
\begin{Lem}\label{Lem.III.1}
Let $f$ be in ${\cal E}(\R)$, let $c>0$, $a>0$ and let
$S=\R\times]0,c[$. Then there exists at most one $k$-temperature
${\cal U}$ on $S$ which is continuous on $\overline{S}$ and satisfies the
conditions that ${\cal U}(x,0)=f(x)$, $x\in\R$ and
$$\dint_0^c\left[\dint_{\R}|{\cal U}(x,t)|e^{-ax^2}|x|^{2k}dx\right]dt<\infty.$$
\end{Lem}
{\footnotesize\bf{Proof}}\quad Since $V_k$ is a topological automorphism to the space ${\cal E}(\R)$, then from Theorem 16 of Friedman \cite{Friedman}
(see also Lemma 5 of Flett \cite{Flett1}), there exists at most one classical temperature $\tilde{{\cal U}}$ on $S$ which is continuous on $\overline{S}$
and satisfies the conditions that
$$\tilde{{\cal U}}(x,0)=V_k^{-1}(f)(x),\,\,\,x\in\R\,\,\,\mbox{and}\,\,\,\dint_0^c\left[\dint_{\R}|\tilde{{\cal U}}(x,t)|e^{-ax^2}dx\right]dt<\infty.$$
Thus, $(x,t)\longmapsto {\cal U}(x,t)=V_k(\tilde{{\cal U}}(.,t))(x)$ is a $k$-temperature on $S$ which is continuous on $\overline{S}$ and ${\cal U}(x,0)=f(x)$, $x\in\R$.
From the formula (\ref{e:I.1}) we deduce that for $x\neq0$
\begin{equation}\label{e.VI.150}
V_k(\tilde{{\cal U}}(.,t))(x)=B(k)|x|^{-2k}sgn(x)\dint_{-|x|}^{|x|}\tilde{{\cal U}}(y,t)(x^2-y^2)^{k-1}(x+y)dy.
\end{equation}
Then according to Fubini-Tonelli's theorem, formula (\ref{e.VI.150}), change of variables $\xi=x^2$ and formula (11) given in \cite{Bateman} p. 202, we have
$$\dint_0^c\left[\dint_{\R}|{\cal U}(x,t)|e^{-ax^2}|x|^{2k}dx\right]dt\leq\dint_0^c\left[\dint_{\R}V_k(|\tilde{{\cal U}}(.,t)|)(x)e^{-ax^2}|x|^{2k}dx\right]dt$$
$$\leq B(k)\dint_0^c\left[\dint_{\R}|\tilde{{\cal U}}(y,t)|\left(\dint_{y^2}^{+\infty}e^{-a\xi}(\xi-y^2)^{k-1}d\xi\right)dy\right]dt$$
$$=B(k,a)\dint_0^c\left[\dint_{\R}|\tilde{{\cal U}}(y,t)|e^{-ay^2}dy\right]dt<\infty.$$ This achieves the proof.
\begin{Th}\label{Th.III.2}
Let $p\in[1,\infty]$ and let ${\cal U}$ be a $k$-temperature on $\R^2_+$
such that the function $t\longmapsto\|{\cal U}(.,t)\|_{k,p}$ is locally
integrable on $]0,\infty[$. Hence
\\ (i) for all $s>0$ and $(x,t)\in\R^2_+$,
\begin{equation}\label{semigroup}
{\cal U}(x,s+t)=\dint_{\R}{\cal T}^k_{-y}F^k_t(x){\cal U}(y,s)|y|^{2k}dy.
\end{equation}
(ii) $t\longmapsto\|{\cal U}(.,t)\|_{k,p}$ is decreasing and
continuous on $]0,\infty[$. Further, for each $(n,m)\in\N\times\N$
the function $t\longmapsto\|{\cal
D}_k^n\partial^m_t{\cal U}(.,t)\|_{k,p}$ is decreasing and continuous on $]0,\infty[$.
\end{Th}
{\footnotesize\bf{Proof}}\quad It is obtained in the same way as for Theorem 4 of Flett \cite{Flett1} by
using Lemma \ref{Lem.III.1}.
\begin{Rem}
The equation (\ref{semigroup}) is called the "semi-group formula"
hereafter.
\end{Rem}
\section{\footnotesize{{\bf{Dunkl-Bessel Potentials}}}}\par The
aim of this section is to define the Bessel potential of some
classes of $k$-temperature associated with the Dunkl setting on $\R$ and to prove related properties needed
later. We adopt the method used by Flett \cite{Flett1} and Johnson \cite{Johnson} in treating classical temperatures.
\begin{Def}
For any $f\in L^p(\R,|x|^{2k}dx)$, where $1\leq p\leq\infty$
and for any $\alpha>0$, the Dunkl-Bessel potential ${\cal
J}^k_{\alpha}f$ of order $\alpha$ of $f$ is given by
$${\cal J}^k_{\alpha}f:={\cal B}^k_{\alpha}\ast_kf,$$ with the kernel
function
\begin{equation}\label{e.IV.1}
\begin{array}{rcl} {\cal B}^k_{\alpha}(x) &:=& \dfrac{1}
{2^{k+\frac{1}{2}}\Gamma(\frac{\alpha}{2})}\dint_0^{+\infty}
e^{-t}e^{-\frac{x^2}{4t}}t^{-k+\frac{(\alpha-1)}{2}-1}dt\\
&=& \dfrac{1}
{2^{\frac{\alpha}{2}-1}\Gamma(\frac{\alpha}{2})}|x|^{
\frac{1}{2}(\alpha-1)-k} K_{\frac{\alpha}{2}-\frac{1}{2}-k}(|x|).
\end{array}
\end{equation}
\end{Def}
   Here $$K_{\beta}(z):=\dfrac{\pi}{2}\left\{\dfrac{J_{-\beta}(z)-J_{\beta}(z)}{\sin\beta\pi}
   \right\},$$ where $J_{\beta}$ is the modified Bessel function
   of the first kind with series expansion
   $$J_{\beta}(z):=\dsum_{n=0}^{+\infty}\dfrac{(\frac{1}{2}z)^{\beta+2n}}{n!\Gamma(\beta+n+1)}.$$
 The Bessel potentials associated with the Dunkl setting on $\R$ which we recall the $k$-Bessel potentials are bounded operators from
   $L^p(\R,|x|^{2k}dx)$ to itself for $1\leq p\leq\infty$ (see
   \cite{Xu2}), i.e., if $f\in L^p(\R,|x|^{2k}dx)$ and $\alpha>0$, then ${\cal J}^k_{\alpha}f\in
   L^p(\R,|x|^{2k}dx)$ and $\|{\cal
   J}^k_{\alpha}f\|_{k,p}\leq\|f\|_{k,p}$. Further, for
   $\alpha,\beta>0$ $${\cal J}^k_{\alpha}({\cal J}^k_{\beta}f)={\cal
   J}^k_{\alpha+\beta}f.$$
   By using the well-known asymptotic behavior of the function $K_{\nu}$, $\nu\in\R$ ( see \cite{Aron} page 415 ), we deduce that
   \\ (a)
   ${\cal B}^k_{\alpha}(x)\sim\dfrac{\Gamma(\frac{1-\alpha}{2}+k)}{2^{\alpha-\frac{1}{2}-k}\Gamma(\frac{\alpha}{2})}|x|^{\alpha-1-2k}$,
   \,\,\footnote{As usual, we write $f(x)\sim g(x)$ as $x\longrightarrow a$ if $\displaystyle\lim_{x\longrightarrow
   a}\frac{f(x)}{g(x)}=1$.}
   as $|x|\longrightarrow0$, for $0<\alpha<2k+1$.
   \\ (b) ${\cal B}^k_{1+2k}(x)\sim\dfrac{1}{2^{k-\frac{1}{2}}\Gamma(k+\frac{1}{2})}\log(\dfrac{1}{|x|})$
   as $|x|\longrightarrow0$.
   \\ (c)
   ${\cal B}^k_{\alpha}(x)\sim\frac{\Gamma(\frac{\alpha-1}{2}-k)}{2^{\frac{1}{2}+k}\Gamma(\frac{\alpha}{2})}$
   as $|x|\longrightarrow0$, for $\alpha>2k+1$.
   \\ (d)
   ${\cal B}^k_{\alpha}(x)\sim\dfrac{\sqrt{\pi}}{2^{\frac{\alpha-1}{2}}\Gamma(\frac{\alpha}{2})}|x|^{\frac{\alpha}{2}-1-k}e^{-|x|}$
   as $|x|\longrightarrow\infty$, for $\alpha>0$.
   \\ As a consequence, we obtain
   \begin{equation}\label{e.IV.2}
   {\cal B}^k_{\alpha}(x)\leq
   B(k,\alpha)|x|^{\alpha-1-2k},\,\,\,\mbox{if}\,\,\,0<\alpha<1+2k.
   \end{equation}
   By differentiation under the integration sign of formula (\ref{e.IV.1}), and using the identity
   $$t^{-a}=\dfrac{1}{\Gamma(a)}\dint_0^{+\infty}e^{-t\delta}\delta^a\dfrac{d\delta}{\delta},\,\,\,\mbox{with}\,\,\,a>0,$$
   we show that
   \begin{equation}\label{e.IV.3}
   |{\cal D}_k{\cal B}^k_{\alpha}(x)|<B(k,\alpha)|x|^{\alpha-2-2k},\,\,\mbox{if}\,\,\,0<\alpha<2k+3.
   \end{equation}
   Added to this, we can see that the kernel ${\cal B}_{\alpha}^k$, $\alpha>0$, satisfies
   \\ (i) ${\cal B}_{\alpha}^k(x)\geq0$, for all $x\in\R$.
   \\ (ii) $\|{\cal B}^k_{\alpha}\|_{k,1}=1$.
   \\ (iii) ${\cal
   F}_k({\cal B}^k_{\alpha})(x)=(1+x^2)^{-\frac{\alpha}{2}}$, $x\in\R$.
   \\ (iv)
   ${\cal B}^k_{\alpha_1+\alpha_2}={\cal B}^k_{\alpha_1}\ast_k{\cal B}^k_{\alpha_2}$,
   if $\alpha_1$, $\alpha_2>0$.

   \par The next theorem is the basis of our definition of the Dunkl-Bessel potential for $k$-temperatures.
   \begin{Th}\cite{N-A-S}\label{Th.IV.1}
   Let $\alpha>0$, $1\leq p\leq\infty$ and let $f\in
   L^p(\R,|x|^{2k}dx)$, then \\ (i) The $k$-Bessel
   potential ${\cal J}^k_{\alpha}f$ of order $\alpha$ of $f$ is
   given for almost all $x$ by
   \begin{equation}\label{e.IV.4}
   {\cal
   J}^k_{\alpha}f(x)=\dfrac{1}
   {\Gamma(\frac{\alpha}{2})}\dint_0^{+\infty}t^{\frac{\alpha}{2}-1}
   e^{-t}G_t^k(f)(x)dt,
   \end{equation}
   where $G_t^k(f)$, $t>0$, is the $k$-heat transform of $f$ on $\R$.
   \\ (ii) The $k$-heat transform of ${\cal
   J}^k_{\alpha}f$, $\alpha>0$,
   on $\R$ is the function $G_s^k({\cal
   J}^k_{\alpha}f)$ given by
   \begin{equation}\label{e.IV.5}
   G_s^k({\cal
   J}^k_{\alpha}f)(x)=\dfrac{1}
   {\Gamma(\frac{\alpha}{2})}\dint_0^{+\infty}t^{\frac{\alpha}{2}-1}
   e^{-t}G^k_{s+t}(f)(x)dt.
   \end{equation}
   Moreover, for each $s>0$, the function $x\mapsto G_s^k({\cal
   J}^k_{\alpha}f)(x)$ is the $k$-Bessel potential of $x\mapsto G_s^k(f)(x)$.
   \end{Th}
   \begin{Def}
    Let ${\cal T}^k(\R^2_+)$ denotes the linear space of $k$-temperatures
    ${\cal U}$ on $\R^2_+$ with the properties that if
    $(n,m)\in\N\times\N$, $b>0$, $c>0$, and $S$ is a compact
    subset of $\R$, then there is a positive constant $C$ such that
    $$
    |{\cal D}_k^{n}\partial^m_t{\cal U}(x,t)|\leq
    Ct^{-b}e^t,\,\,\,\mbox{for all}\,\,(x,t)\in
    S\times[c,\infty[.
    $$
   \end{Def}
   \begin{Def}\label{Def.IV.2}
   For any ${\cal U}$ in ${\cal T}^k(\R^2_+)$ and any real number
   $\alpha$, ${\cal J}^k_{\alpha}{\cal U}$ is the function defined on
   $\R^2_+$ by
   \\ (i) ${\cal
    J}^k_0({\cal U})={\cal U}$;\\ (ii) if $\alpha>0$,
    $$
    {\cal
    J}^k_{\alpha}({\cal U})(x,s)=\dfrac{1}
    {\Gamma(\frac{\alpha}{2})}\dint_0^{+\infty}t^{\frac{\alpha}{2}
    -1}e^{-t}{\cal U}(x,s+t)dt;
    $$
    \\ (iii) if $\alpha$ is a negative even integer, say
    $\alpha=-2m$, then $${\cal J}^k_{\alpha}({\cal U})(x,s)={\cal J}^k_{-2m}({\cal U})(x,s)=(-1)^me^s\partial
    ^m_s\{e^{-s}{\cal U}(x,s)\};$$
    (iv) if $\alpha=-\beta<0$ and $\beta$ is not an even
    integer, then $${\cal J}^k_{\alpha}({\cal U})={\cal J}^k_{-\beta}({\cal U})={\cal J}^k_{2m-\beta}
    \left({\cal J}^k_{-2m}({\cal U})\right);$$ where
    $m=[\frac{1}{2}\beta]+1$, \footnote{Here $[x]$ stands for the greatest integer not
    exceeding $x$, $x\in\R$. } and where ${\cal J}^k_{2m-\beta}$
    and ${\cal J}^k_{-2m}$ are defined as in (ii) and
    (iii).
   \end{Def}
   \begin{Th}\label{Th.IV.2}\cite{N-A-S}
    Let ${\cal U}\in{\cal T}^k(\R^2_+)$ and $\alpha$, $\beta$ be real numbers. \\ (i)
     ${\cal J}^k_{\alpha}({\cal U})$ is well-defined and ${\cal J}^k_{\alpha}({\cal U})\in{\cal
    T}^k(\R^2_+)$, \\ (ii) ${\cal J}^k_{\alpha}
    \left({\cal J}^k_{\beta}({\cal U})\right)={\cal
    J}^k_{\alpha+\beta}({\cal U})={\cal J}^k_{\beta}
    \left({\cal J}^k_{\alpha}({\cal U})\right)$.
    \end{Th}
    \begin{Cor}\label{Cor.IV.01}
    For each real number $\alpha$, ${\cal J}^k_{\alpha}$ is a linear isomorphism of ${\cal T}^k(\R^2_+)$ onto itself, with inverse ${\cal J}^k_{-\alpha}$.
    \end{Cor}
    \begin{Th}\label{Th.IV.3}
    Let $f$ be in $L^p(\R,|x|^{2k}dx)$, $1\leq p\leq\infty$,
    $\alpha>0$, and let $G^k_t(f)$ be the $k$-heat transform of
    $f$ on $\R^2_+$. Then for $t>0$
    \\ (i) $\|{\cal J}^k_{\alpha}G^k_t(f)\|_{k,p}\leq
    c_k^{-1}\|f\|_{k,p}$;
    \\ (ii) $\|{\cal J}^k_{-\alpha}G^k_t(f)\|_{k,p}\leq
    B(k,\alpha)(t^{-\frac{1}{2}\alpha}+1)\|f\|_{k,p}$;
    \\ (iii) furthermore, if $1\leq p<\infty$ then
    $$\|{\cal J}^k_{-\alpha}G^k_t(f)\|_{k,p}=\circ(t^{-\frac{1}{2}\alpha}),\,\,\,\mbox{as}\,\,\,t\longrightarrow0^+.$$
    \end{Th}
  {\footnotesize\bf{Proof}}\quad Part (i) follows from relation
  (\ref{e.IV.5}), Minkowski's integral inequality and Theorem
  \ref{Th.II.2}(i). According to the fact that
  $$J^k_{-2m}G^k_t(f)=\dsum_{i=0}^m(-1)^i(\begin{array}[c]{c}m\\i\end{array})\partial^i_tG^k_t(f),\,\,\,m\in\N,$$
  Minkowski's inequality, Theorem \ref{Th.II.2}(i) and the following inequality
  \begin{equation}\label{Relation}
  (a+b)^s\leq2^{s-1}(a^s+b^s),\,\,s\in[1,+\infty[,\,\,a,b\geq0,
  \end{equation}
   yield the part (ii) when $\alpha=2m$. Supposing that
  $\alpha$ is not an even integer and let
  $m=[\frac{1}{2}\alpha]+1$. Then for $(x,s)\in\R^2_+$
  $${\cal J}^k_{-\alpha}G^k_s(f)(x)=\dfrac{1}{\Gamma(m-\frac{1}{2}\alpha)}\dint_0^{+\infty}t^{m-\frac{1}{2}\alpha-1}e^{-t}
  {\cal J}^k_{-2m}G^k_{s+t}(f)(x)dt.$$ Hence, Minkowski's integral
  inequality and the previous case when $\alpha=2m$ yield that $\|{\cal J}^k_{-\alpha}G^k_s(f)\|_{k,p}\leq
  B(k,\alpha)(s^{-\frac{1}{2}\alpha}+1)\|f\|_{k,p}$. We shall prove (iii)
  only when $\alpha=2m$, because the general case can be treated
  in the same manner. Let $(x,t)$ be in $\R^2_+$. Thus by
  Proposition \ref{Prop.II.1}(iv)
  $${\cal J}^k_{-2m}G^k_t(f)(x)=\dsum_{i=0}^m(-1)^i(\begin{array}[c]{c}m\\i\end{array})\dint_{\R}\partial^i_tF^k_t(y)({\cal
  T}^k_{-y}f(x)-f(x))|y|^{2k}dy$$ which together with Minkowski's integral
  inequality imply that
  $$t^m\|{\cal J}^k_{-2m}G^k_t(f)\|_{k,p}\leq t^m\dsum_{i=0}^m(\begin{array}[c]{c}m\\i\end{array})
  \dint_{|y|<\delta}|\partial^i_tF^k_t(y)|\|{\cal T}^k_{-y}f-f\|_{k,p}|y|^{2k}dy+$$
  $$t^m\dsum_{i=0}^m(\begin{array}[c]{c}m\\i\end{array})\dint_{|y|\geq\delta}|\partial^i_tF^k_t(y)|\|{\cal
  T}^k_{-y}f-f\|_{k,p}|y|^{2k}dy=I_1(t)+I_2(t)\,\,\,(\delta>0).$$
  Since $\lim_{y\longrightarrow0}\|{\cal T}^k_{-y}f-f\|_{k,p}=0$,
  for an arbitrary positive number $\epsilon$, there exists a
  $\delta>0$ such that $\|{\cal T}^k_{-y}f-f\|_{k,p}<\epsilon$ if
  $|y|<\delta$. Therefore, from Proposition \ref{Prop.II.1}(iv) and inequality (\ref{Relation}), we obtain $I_1(t)\leq B(k,m)(1+t^m)\epsilon$. By
  relation (\ref{e.I.4}), Proposition \ref{Prop.II.1}(iii) and the change of variables, we have
  $$I_2(t)\leq
  B(k)\|f\|_{k,p}\dsum_{i=0}^m(\begin{array}[c]{c}m\\i\end{array})t^{m-i}\dint_{\frac{\delta^2}{4t}}^{+\infty}|R_i(\sigma)|
  e^{-\sigma}\sigma^{k-1/2}d\sigma.$$
  Letting $t\rightarrow0^+$, the last integral approaches to $0$. This proves the part (iii).
  \begin{Cor}\label{Cor.IV.1}
  Let $\alpha>0$, $1\leq p\leq\infty$, and ${\cal U}$ be in ${\cal
  T}^k(\R^2_+)$. If ${\cal U}$ satisfies the semi-group formula, then
  for all $s,t>0$
  \\ (i) $\|{\cal
  J}^k_{\alpha}{\cal U}(.,s+t)\|_{k,p}\leq\|{\cal U}(.,s)\|_{k,p}$.
  \\ (ii) $\|{\cal
  J}^k_{-\alpha}{\cal U}(.,s+t)\|_{k,p}\leq B(k,\alpha)(t^{-\frac{1}{2}\alpha}+1)\|{\cal U}(.,s)\|_{k,p}$.
  \end{Cor}
  {\footnotesize\bf{Proof}}\quad Let $s$ be fixed. We may assume that
  $\|{\cal U}(.,s)\|_{k,p}$ is finite (otherwise the conclusion would be trivial). Then for all $t>0$, by the semi-group formula for ${\cal U}$ yields
  $${\cal U}(x,s+t)=\dint_{\R}{\cal T}^k_{-y}F^k_t(x){\cal U}(y,s)|y|^{2k}dy$$
  which implies the corollary by analogous reasoning of Theorem \ref{Th.IV.3}.
  \begin{Th}\label{Th.IV.4}
  Let $1\leq p\leq\infty$, $1\leq q<\infty$, $\beta$ be a positive
  number and ${\cal U}$ be a $k$-temperature on $\R^2_+$ such that
  $$C=\left\{\dint_0^{+\infty}t^{\frac{1}{2}q\beta-1}e^{-t}\|{\cal U}(.,t)\|^q_{k,p}dt\right\}^{\frac{1}{q}}<\infty.$$
  Thus for $t>0$, $\|{\cal U}(.,t)\|_{k,p}\leq
  B(q,\beta)(1+t^{-\frac{1}{2}\beta})C$ and
  $\|{\cal U}(.,t)\|_{k,p}=\circ(t^{-\frac{1}{2}\beta})$ as $t\longrightarrow0^+$.
  Moreover, if $q<r<\infty$, then
  $$\left\{\dint_0^{+\infty}t^{\frac{1}{2}r\beta-1}e^{-t}\|{\cal U}(.,t)\|^r_{k,p}dt\right\}^{\frac{1}{r}}\leq B(q,r,\beta)C.$$
  \end{Th}
  {\footnotesize\bf{Proof}}\quad The proof is similar to the classical case (see Theorem 11 p. 405 in \cite{Flett1}).
  \begin{Th}\label{Th.IV.5}
  Let $1\leq p\leq\infty$, $1\leq q\leq\infty$, $\alpha$ be a real
  number, $\beta>0$, $\beta>\alpha$ and ${\cal U}$ be a
  $k$-temperature on $\R^2_+$ such that
  $$C:=\left\{
  \begin{array}{rcl}
  \left\{\dint_0^{+\infty}t^{\frac{1}{2}q\beta-1}e^{-t}\|{\cal U}(.,t)\|^q_{k,p}dt\right\}^{\frac{1}{q}}=C_1<\infty,\,\,(1\leq q<\infty),\\
  \displaystyle\sup_{t>0}\left\{t^{\frac{1}{2}\beta}e^{-t}\|{\cal U}(.,t)\|_{k,p}\right\}=C_2<\infty,\,\,(q=\infty).\,\,\,\,\,\,\,\,\,\,\,\,\,\,\,\,\,\,\,\,\,\,\,\,\,
  \end{array}
  \right.
  $$
  (i) ${\cal U}\in{\cal T}^k(\R^2_+)$ and
  $$\left\{
  \begin{array}{rcl}
  \left\{\dint_0^{+\infty}t^{\frac{1}{2}q(\beta-\alpha)-1}e^{-t}\|{\cal J}^k_{\alpha}{\cal U}(.,t)\|^q_{k,p}dt\right\}^{\frac{1}{q}}\leq
  B(k,\alpha,\beta,q)C_1,
  \,\,(1\leq q<\infty),\\ \displaystyle\sup_{t>0}\left\{t^{\frac{1}{2}(\beta-\alpha)}e^{-t}\|{\cal J}^k_{\alpha}{\cal U}(.,t)\|_{k,p}\right\}
  \leq B(k,\alpha,\beta)C_2,\,(q=\infty).\,\,\,\,\,\,
  \end{array}
  \right.
  $$
  (ii) If $1\leq q<\infty$, then $\|{\cal
  J}^k_{\alpha}{\cal U}(.,t)\|_{k,p}=\circ(t^{-\frac{1}{2}(\beta-\alpha)})$
  as $t\longrightarrow0^+$.
 \\ (iii) If $q=\infty$ and
  $\|{\cal U}(.,t)\|_{k,p}=\circ(t^{-\frac{1}{2}\beta})$ as $t\longrightarrow0^+$,
  then $\|{\cal
  J}^k_{\alpha}{\cal U}(.,t)\|_{k,p}=\circ(t^{-\frac{1}{2}(\beta-\alpha)})$
  as $t\longrightarrow0^+$.
  \end{Th}
   {\footnotesize\bf{Proof}}\quad Clearly $t\longmapsto\|{\cal U}(.,t)\|_{k,p}$ is locally
  integrable on $]0,\infty[$, so that ${\cal U}\in{\cal T}^k(\R^2_+)$ and
  $\|{\cal U}(.,t)\|_{k,p}$ is decreasing. Therefore ${\cal
  J}^k_{\alpha}{\cal U}$ is well defined. First, suppose that
  $\gamma=-\alpha>0$. Then by Corollary \ref{Cor.IV.1} we see that
  \begin{equation}\label{e.IV.6}
  \|{\cal J}^k_{\alpha}{\cal U}(.,2t)\|_{k,p}\leq
  B(k,\alpha)(t^{\frac{1}{2}\alpha}+1)\|{\cal U}(.,t)\|_{k,p}
  \end{equation}
  which implies that
  $$\left\{\dint_0^{+\infty}t^{\frac{1}{2}q(\beta-\alpha)-1}e^{-t}\|{\cal J}^k_{\alpha}{\cal U}(.,t)\|^q_{k,p}dt\right\}^{\frac{1}{q}}\leq B(k,\alpha,\beta,q)C_1.$$
  Next, we shall prove the result for the special case when $\alpha=2$ and
  $\beta>2$. Since
 \begin{equation}\label{e.IV.7}
 {\cal J}^k_2{\cal U}(x,t)=\dint_0^{+\infty}e^{-\xi}{\cal U}(x,t+\xi)d\xi,
 \end{equation}
  it follows from Minkowski's integral inequality and Hardy's inequality that
  $$\left\{\dint_0^{+\infty}t^{\frac{1}{2}q(\beta-2)-1}e^{-qt}\|{\cal J}^k_2{\cal U}(.,t)\|^q_{k,p}dt\right\}^{\frac{1}{q}}\leq
  B(k,\beta,q)\left\{\dint_0^{+\infty}t^{\frac{1}{2}q\beta-1}e^{-qt}\|{\cal U}(.,t)\|^q_{k,p}dt\right\}^{1/q}.$$
  To prove the result for $\alpha=\delta>0$, let
  $\gamma$ be the least positive number such that $\gamma+\delta$
  is an even positive integer. Then by applying part (i) in case $\alpha<0$, we
  have
  $$\left\{\dint_0^{+\infty}t^{\frac{1}{2}q(\beta+\gamma)-1}e^{-t}\|{\cal J}^k_{-\gamma}{\cal U}(.,t)\|^q_{k,p}dt\right\}^{\frac{1}{q}}\leq B(k,\gamma,\beta,q)C_1$$
  and hence after repeated applications of part (i) in case $\alpha=2$,
  we obtain
  $$\left\{\dint_0^{+\infty}t^{\frac{1}{2}q(\beta-\delta)-1}e^{-t}\|{\cal J}^k_{\delta} {\cal U}(.,t)\|^q_{k,p}dt\right\}^{\frac{1}{q}}\leq
  B(k,\alpha,\beta,q)C_1.$$ It is easy to see
  $$\displaystyle\sup_{t>0}\left\{t^{\frac{1}{2}(\beta-\alpha)}e^{-t}\|{\cal J}^k_{\alpha}{\cal U}(.,t)\|_{k,p}\right\}\leq B(k,\alpha,\beta)C_2$$ from Corollary \ref{Cor.IV.1}.
  The assertion (ii) then follows from part (i) and Theorem
  \ref{Th.IV.4}.\\ Now, we shall prove the assertion (iii). First, assuming that
  $\alpha<0$, the result follows easily from the estimate
  (\ref{e.IV.6}). Next, we shall prove the result for the case when
  $\alpha=2$ and $\beta>2$. It follows from relation (\ref{e.IV.7}) and Minkowski's integral inequality
  that
  $$s^{\frac{\beta}{2}-1}\|{\cal J}^k_2{\cal U}(.,s)\|_{k,p}\leq
  s^{\frac{\beta}{2}-1}e^s\dint_s^{+\infty}e^{-t}\|{\cal U}(.,t)\|_{k,p}dt,$$
  consequently the assertion is proved for the special case. In case
  $\alpha=\delta>0$ and by choosing $\gamma>0$, $\gamma+\delta$ is
  an even positive integer. Applying the above result for $\alpha<0$ we see
  that $\|{\cal
  J}^k_{-\gamma}{\cal U}(.,t)\|_{k,p}=\circ(t^{-\frac{1}{2}(\beta+\gamma)})$.
  Repeated use of the result for $\alpha=2$ yields
  $\|{\cal J}^k_{\delta}{\cal U}(.,t)\|_{k,p}=\|{\cal J}^k_{\gamma+\delta}
  ({\cal J}^k_{-\gamma}{\cal U}(.,t))\|_{k,p}=\circ(t^{-\frac{1}{2}(\beta+\gamma)+\frac{1}{2}(\gamma+\delta)})
  =\circ(t^{-\frac{1}{2}(\beta-\delta)})$. Thus part (iii) is proved.
  \begin{Def}
For any real number $\alpha$ and for any $T\in S'(\R)$, the
$k$-Bessel potential of order $\alpha$ of $T$ is the element
${\cal J}^k_{\alpha}(T)$ of $S'(\R)$ given by the relation
$${\cal F}_k({\cal J}^k_{\alpha}(T)):=(1+(.)^2)^{-\frac{\alpha}{2}}{\cal F}_k(T),$$
where the identity is to be understood in the sense of
distributions.
\end{Def}
\begin{Rems} We have
\begin{itemize}
\item For all real $\alpha$, $\beta$ and all $T\in S'(\R)$
$${\cal J}^k_{\alpha}({\cal J}^k_{\beta}(T))={\cal J}^k_{\alpha+\beta}(T).$$
\item By definition
$${\cal J}^k_{\alpha}(T)=T\ast_k{\cal B}^k_{\alpha},$$
where ${\cal B}^k_{\alpha}$ is a tempered distribution whose Dunkl
transform ${\cal
F}_k({\cal B}^k_{\alpha})=\left[(1+(.)^2)^{-\frac{\alpha}{2}}\right]$,
\footnote{$[f]$ is the distribution on $\R$ associated with the
function $f$. In addition $[f]$ belongs to $S'(\R)$, when $f\in
L^p(\R,|x|^{2k}dx)$ or $f$ is slowly increasing}.
\item If $f\in L^p(\R,|x|^{2k}dx)$, where $p\in [1,\infty]$ and
$\alpha>0$, then
$${\cal J}^k_{\alpha}([f])={\cal J}^k_{\alpha}(f)=f\ast_k{\cal B}^k_{\alpha}.$$
\end{itemize}
\end{Rems}

  \section{\footnotesize{{\bf{Generalized Dunkl-Lipschitz Spaces, $\alpha$ Real}}}}
\par Our basic aim is to define Lipschitz spaces associated with the Dunkl
operators for all real $\alpha$. In the classical case, the heat (or
Poisson) semi-group provides an alternative characterization of
the Lipschitz spaces, we will follow this approach, using the
$k$-heat (or $k$-Poisson) semi-group, to define generalized Dunkl-Lipschitz
spaces. One of the main result of this part is to show that ${\cal J}^k_{\beta}$ is an isomorphism of $\wedge^k_{\alpha,p,q}(\R)$ onto $\wedge^k_{\alpha+\beta,p,q}(\R)$ for real $\alpha$ and $\beta$. The section closes by giving some properties and continuous embedding for the space $\wedge^k_{\alpha,p,q}(\R)$.
\par We define for $t>0$, the function $P^k_t$ on $\R$ by
  $$P^k_t(x):=\tilde{c}_k\,\dfrac{t}{(t^2+x^2)^{k+1}},\,\,\,\mbox{where}\,\,\,\,\tilde{c}_k :=\dfrac{2^{k+\frac{1}{2}}}{\Gamma(\frac{1}{2})}\Gamma(k+1).$$
  The function $P^k_t$ is called the $k$-Poisson kernel. We summarize the properties of $P^k_t$ in the following proportion :
  \begin{Prop}\label{Prop.V.1}
  For all $t>0$, $n\in\N$ and $x\in\R$, we have
  \\ (i) ${\cal F}_k(P^k_t)(x)=e^{-t|x|}$.
  \\ (ii) $\int_{\R}P^k_t(y)|y|^{2k}dy=1$.
  \\ (iii) $P^k_t\in L^p(\R,|x|^{2k}dx)$, $1\leq p\leq\infty$.
  \\ (iv) $P^k_{t_1+t_2}=P^k_{t_1}\ast_kP^k_{t_2}$, if $t_1,t_2>0$.
  \\ (v) $\|\partial^n_tP^k_t\|_{k,1}\leq B(k,n)t^{-n}$,  $\|{\cal D}^n_kP^k_t\|_{k,1}\leq \tilde{B}(k,n)t^{-n}$ and $|\partial^n_tP^k_t(x)|\leq B(k,n)
  t^{-2k-1-n}$.
  \\ (vi) $\displaystyle\lim_{t\rightarrow0}P^k_tf(x)=f(x)$, where the limit is interpreted in $L^k_p$-norm and pointwise a.e. For $f\in C_0(\R)$ the convergence is uniform on $\R$.
  \end{Prop}
  However, for $t>0$ and for all $f\in L^p(\R,|x|^{2k}dx)$, $p\in[1,\infty]$, we put
  $$P^k_tf(x):=P^k_t\ast_kf(x),\,\,x\in\R.$$
  The function $P^k_tf$ is called the Poisson transform of a function
  $f$ associated with the Dunkl setting on $\R$ that's why we may recall it the $k$-Poisson transform of $f$.
  \\ A $C^2$ function ${\cal U}$ on $\R^2_+$ satisfying $({\cal D}^2_k+\partial^2_t){\cal U}(x,t)=0$ is said to be $k$-harmonic. For $p\in[1,\infty]$, we suppose that
  \begin{equation}\label{e.VI.100}
  A^p:=\displaystyle\sup_{t>0}B(k)\dint_{\R}|{\cal U}(x,t)|^p|x|^{2k}dx<\infty.
  \end{equation}
  Now, we need the following key results.
  \begin{Lem} (Semi-group property)\label{Lem.VI.Semi-group}
  If ${\cal U}(x,t)$ is $k$-harmonic on $\R^2_+$ and bounded in each proper sub-half space of $\R^2_+$, then for $t_0>0$, ${\cal U}(x,t+t_0)$ is identical with the $k$-Poisson transform of ${\cal U}(.,t_0)$, that is,
  $${\cal U}(x,t_0+t)=P^k_t({\cal U}(.,t_0))(x),\,\,\mbox{for}\,\,t>0.$$
  Furthermore,
  $$\partial_t{\cal U}(x,t_0+t)=\partial_tP^k_t({\cal U}(.,t_0))(x)=P^k_t(\partial_t{\cal U}(.,t_0))(x).$$
  \end{Lem}
  {\footnotesize\bf{Proof}}\quad It is obtained in the same way as for property 12 p. 417 in \cite{Taibleson}.
  \begin{Th} \label{Th.VI.characterization}(Characterization of $k$-Poisson transform)
  Let $p\in[1,\infty]$ and let ${\cal U}(x,t)$ be $k$-harmonic on $\R^2_+$. Then
  \\ (i) if $1<p<\infty$, ${\cal U}(x,t)$ is the $k$-Poisson transform of a function $f\in L^p(\R,|x|^{2k}dx)$ if and only if ${\cal U}(x,t)$ satisfies condition (\ref{e.VI.100}), moreover $\|f\|_{k,p}=A$.
  \\ (ii) For $p=1$, ${\cal U}(x,t)$ is the $k$-Poisson transform of $f\in L^1(\R,|x|^{2k}dx)$ if and only if ${\cal U}(x,t)$ satisfies condition (\ref{e.VI.100}) and
  $\|{\cal U}(.,t_1)-{\cal U}(.,t_2)\|_{k,1}$, as $t_1,t_2\rightarrow0$.
  \\ (iii) For $p=\infty$, ${\cal U}(x,t)$ is the $k$-Poisson transform of a function $f\in L^{\infty}(\R,|x|^{2k}dx)$ if and only if there exists $C>0$ such that
  $\|{\cal U}(.,t)\|_{k,\infty}\leq C$ for all $t>0$.
  \end{Th}
  {\footnotesize\bf{Proof}}\quad Parts (i) and (ii) are proved in \cite{Zh-Ji} Theorem 4.16 p. 254. Part (iii) is proved in usual way (see \cite{Taibleson} p. 416).
  \begin{Rem}\label{Rem.VI.Temperature}
  Analogously to the $k$-harmonic case, we can assert that Theorem \ref{Th.VI.characterization} and Lemma \ref{Lem.VI.Semi-group} are true when we take ${\cal U}(x,t)$
  $k$-temperature on $\R^2_+$ and we replace $k$-Poisson transform by $k$-heat transform.
  \end{Rem}
\par Before giving a central result of this section, we need to
recall the definition of the spaces $\wedge^k_{\alpha,p,q}(\R)$ (see \cite{S-K}) and the following auxiliary lemmas.
\begin{Def}
  The generalized Dunkl-Lipschitz spaces $\wedge^k_{\alpha,p,q}(\R)$,
  $\alpha\in]0,1[$, $1\leq p,q\leq\infty$, is the set of functions $f\in
  L^p(\R,|x|^{2k}dx)$ for which the norm
  $$\|f\|_{k,p}+\left\{\dint_{\R}\dfrac{\|\triangle_{y,k}f\|^q_{k,p}}{|y|^{1+\alpha
  q}}dy\right\}^{\frac{1}{q}}<\infty,\footnote{$\triangle_{y,k}f={\cal T}^k_yf-f$}\,\,\mbox{if}\,\,q<\infty$$
  and $$\|f\|_{k,p}+\displaystyle\sup_{|y|>0}\dfrac{\|\triangle_{y,k}f\|_{k,p}}{|y|^{\alpha}}<\infty,\,\,\mbox{if}\,\,q=\infty.$$
  \end{Def}
  {\footnotesize{\bf{Notations}}}
   \begin{itemize}
  \item For any $k$-harmonic (or $k$-temperature) ${\cal U}$ on $\R^2_+$, we denote by
   $${\cal A}^k_{p,q}({\cal U}):=\left\{
  \begin{array}{rcl}
  \left\{\dint_0^{\infty}\left[\|{\cal
  U}(.,t)\|_{k,p}\right]^q\dfrac{dt}{t}\right\}^{\frac{1}{q}}&(1\leq
  q<\infty),\\ \displaystyle\sup_{t>0}\|{\cal
U}(.,t)\|_{k,p}\,\,\,\,\,\,\,\,\,\,\,\,\,\,\,\,&(q=\infty),
\end{array}
 \right.$$
 and
 $${\cal A}^{k,\ast}_{p,q}({\cal U}):=\left\{
\begin{array}{rcl}
\left\{\dint_0^1\left[\|{\cal
U}(.,t)\|_{k,p}\right]^q\dfrac{dt}{t}\right\}^{\frac{1}{q}}&(1\leq
q<\infty),\\ \displaystyle\sup_{0<t\leq1}\|{\cal
U}(.,t)\|_{k,p}\,\,\,\,\,\,\,\,\,\,\,\,\,\,\,\,&(q=\infty),
\end{array}
 \right.$$ the value $\infty$ being allowed.
\item For $\alpha$ real,
$\overline{\alpha}$ will denote the smallest non-negative integer
larger than $\alpha$.
\end{itemize}
\begin{Rems}(\cite{S-K})  We have :
  \begin{itemize}
  \item For $\alpha\in]0,1[$ and $q=\infty$, $f\in\wedge^k_{\alpha,p,\infty}(\R)$ if and only if $\|\partial_tP^k_tf\|_{k,p}
  \leq B(k,\alpha)t^{-1+\alpha}.$

  \item For $\alpha>0$, $p,q\in[1,\infty]$, we set
  $$\wedge^k_{\alpha,p,q}(\R):=\left\{f\in L^p(\R,|x|^{2k}dx):\,\,
  {\cal A}^k_{p,q}(t^{\overline{\alpha}-\alpha}\partial^{\overline{\alpha}}_tP^k_t(f))<\infty\right\}.$$
  The $\wedge^k_{\alpha,p,q}$-norms are defined by
  $$\|f\|_{\wedge^k_{\alpha,p,q}}:=\|f\|_{k,p}+{\cal A}^k_{p,q}(t^{\overline{\alpha}-\alpha}\partial^{\overline{\alpha}}_tP^k_t(f))
  .$$
  \end{itemize}
  \end{Rems}
\begin{Lem}\label{Lem.V.1}
We have
$${\cal B}^k_{\alpha}\in\wedge^k_{\alpha,1,\infty}(\R),\,\,\,\mbox{if}\,\,\,\alpha>0.$$
\end{Lem}
{\footnotesize\bf{Proof}}\quad Let us first consider the case
$\alpha\in]0,1[$. Since ${\cal B}^k_{\alpha}\in L^1(\R,\,|x|^{2k}dx)$, we
can write
$$\|\triangle_{y,k}{\cal B}^k_{\alpha}\|_{k,1}=
\dint_{|x|\leq 2|y|}|{\cal T}^k_y{\cal B}^k_{\alpha}(x)-{\cal B}^k_{\alpha}(x)||x|^{2k}dx+ \dint_{|x|>
2|y|}|{\cal T}^k_y{\cal B}^k_{\alpha}(x)-{\cal B}^k_{\alpha}(x)||x|^{2k}dx=I_1(y)+I_2(y).$$
 ${\cal B}^k_{\alpha}$ is an even function, then formula (\ref{e.I.3})
yields
$${\cal T}^k_y{\cal B}^k_{\alpha}(x)=d_k\dint_0^{\pi}{\cal B}^k_{\alpha}(G(x,y,\theta))h^e(x,y,\theta)\sin^{2k-1}\theta d\theta$$
which shows that ${\cal T}^k_y{\cal B}^k_{\alpha}(x)\geq0$ since ${\cal B}^k_{\alpha}$ is non-negative.
Moreover, using the following inequalities $G(x,y,\theta)\geq||x|-|y||$, $0\leq h^e(x,y,\theta)\leq2$ and relation
(\ref{e.IV.1}), we have
\begin{equation}\label{e.V.2}
{\cal T}^k_{y}{\cal B}^k_{\alpha}(x)\leq 2{\cal B}^k_{\alpha}(|x|-|y|).
\end{equation}
Then, by inequalities (\ref{e.V.2}) and (\ref{e.IV.2}), we have
$$I_1(y)\leq B(k,\alpha)\left\{\dint_{|x|\leq2|y|}||x|-|y|||^{\alpha-1-2k}|x|^{2k}dx+\dint_{|x|\leq2|y|}|x|^{\alpha-1}dx\right\}\leq B(k,\alpha)|y|^{\alpha}.$$
By the generalized Taylor formula with integral remainder
(\ref{e.II.Taylor}), we have
\begin{equation}\label{e.V.12'}
|{\cal
T}^k_y{\cal B}^k_{\alpha}(x)-{\cal B}^k_{\alpha}(x)|\leq\dint_{-|y|}^{|y|}|{\cal
T}^k_z({\cal D}_k{\cal B}^k_{\alpha})(x)|dz.
\end{equation}
 Since ${\cal
D}_k{\cal B}^k_{\alpha}$ is an odd function, formula (\ref{e.I.3}) gives
$${\cal T}^k_z({\cal D}_k{\cal B}^k_{\alpha})(x)=d_k\dint_0^{\pi}{\cal D}_k{\cal B}^k_{\alpha}(G(x,z,\theta))h^o(x,z,\theta)\sin^{2k-1}\theta d\theta.$$
It is obvious to see that $h^o(x,z,\theta)\leq2$ and $0\leq
G(x,z,\theta)\leq|x|+|z|$. Thus, formula (\ref{e.IV.3}) yields
$$|{\cal T}^k_{z}({\cal D}_k{\cal B}^k_{\alpha})(x)|\leq B(k,\alpha)(|x|+|z|)^{\alpha-2-2k}.$$
Hence, by relation (\ref{e.V.12'}) we obtain
$$|{\cal T}^k_y{\cal B}^k_{\alpha}(x)-{\cal B}^k_{\alpha}(x)|\leq B(k,\alpha)|y||x|^{\alpha-2-2k}$$
and so $I_2(y)\leq B(k,\alpha)|y|^{\alpha}$. This completes the proof  when $\alpha\in]0,1[$. To pass to
the general case for $\alpha>0$, we write $t=t_1+t_2+\cdots +t_{\overline{\alpha}}$ and $t_i>0$. Then
$$P^k_t{\cal B}^k_{\alpha}=P^k_{t_1}{\cal B}^k_{\beta}\ast_kP^k_{t_2}{\cal B}^k_{\beta}\ast_k\cdots\ast_kP^k_{t_{\overline{\alpha}}}{\cal B}^k_{\beta},$$
where $\beta=\frac{\alpha}{\overline{\alpha}}\in]0,1[$. Therefore $\|\partial^{\overline{\alpha}}_tP^k_t{\cal B}^k_{\alpha}\|_{k,1}\leq
B(k,\alpha)t^{\alpha-\overline{\alpha}}$, whenever
$t_1=t_2=\cdots=t_{\overline{\alpha}}=\frac{t}{\overline{\alpha}}$. This finishes the proof.
\begin{Lem}\label{Lem.V.2}
Let $1\leq p,q\leq\infty$,  ${\cal U}(x,t)$ is  $k$-harmonic on $\R^2_+$ and bounded in each
proper sub-half space of $\R^2_+$. Suppose we are given $A>0$,
$\alpha>0$, $t_0>0$ and an integer $n>\alpha$ such that
$${\cal A}^k_{p,q}(t^{n-\alpha}\partial^n_t{\cal U})\leq A,$$$$\|{\cal U}(.,t)\|_{k,p}\leq A,\,\,\,t\geq t_0.$$
Then ${\cal U}(x,t)$ is the $k$-Poisson transform of a function
$f\in\wedge^k_{\alpha,p,q}(\R)$ and :
\\ (a) $\|\partial_t{\cal U}(.,t)\|_{k,p}=o(t^{-1})$,
as $t\longrightarrow0$,
\\ (b) $\|f\|_{\wedge^k_{\alpha,p,q}}\leq
B(\alpha,k,t_0,n)A$.
\end{Lem}
{\footnotesize\bf{Proof}}\quad Consider
first the case $\alpha\in]0,1[$. We are given ${\cal U}(.,t)=O(1)$, \footnote{$f(x)=O(g(x))$, $x\rightarrow a$, means $\frac{f(x)}{g(x)}$ is
bounded as $x\rightarrow a$.} as
$t\longrightarrow\infty$, then from Lemma \ref{Lem.VI.Semi-group}, H\"older's inequality and Proposition \ref{Prop.V.1}(v), we get $\partial^{m-1}_t{\cal U}(.,t)=\circ(1)$, $t\longrightarrow\infty$, $m\in\N$.
Using the fact that
$$\partial^{m-1}_t{\cal U}(x,t)=-\dint_t^{\infty}\partial^m_s{\cal U}(x,s)ds,\,\,m\in\N,$$ and
Minkowski's integral inequality, we obtain
\begin{equation}\label{e.V.2'}
\|\partial^{m-1}_t{\cal U}(.,t)\|_{k,p}\leq\dint_t^{+\infty}\|\partial^m_s{\cal U}(.,s)\|_{k,p}ds.
\end{equation}
From Hardy inequality and relation (\ref{e.V.2'}), we
deduce that
$${\cal A}^k_{p,q}(t^{1-\alpha}\partial_t{\cal U})\leq
B(n,\alpha)A.$$ But $t\longmapsto\|\partial_t{\cal U}(.,t)\|_{k,p}$ is a
non-increasing function, so that
$$((1-\alpha)q)^{-\frac{1}{q}}s^{1-\alpha}\|\partial_s{\cal U}(.,s)\|_{k,p}=\left[\dint_0^s(t^{1-\alpha}\|\partial_s{\cal U}(.,s)\|_{k,p})^q
\dfrac{dt}{t}\right]^{\frac{1}{q}}\leq
B(n,\alpha)A,\,\,\mbox{if}\,\,\alpha<1,$$ which proves
\begin{equation}\label{e.V.2"}
t\|\partial_t{\cal U}(.,t)\|_{k,p}\leq
B(n,\alpha,q)At^{\alpha}=\circ(1),\,\,\mbox{as}\,\,t\longrightarrow0.
\end{equation}
If $\alpha\geq1$, it is easily to verify that
\begin{equation}\label{e.V.100}
{\cal A}^k_{p,q}(t^{n-\frac{1}{2}}\partial^n_t{\cal U})\leq B(n,k,q,t_0)A.
\end{equation}
Then by relation (\ref{e.V.2'}), Hardy inequality and relation (\ref{e.V.100}), we have
$${\cal A}^k_{p,q}(t^{\frac{1}{2}}\partial_t{\cal U})\leq B(n,k,q,t_0)A.$$
By the same reason for
$\alpha\in]0,1[$, we obtain $t\|\partial_t{\cal U}(.,t)\|_{k,p}=\circ(1)$ as
$t\longrightarrow0^+$ which proves the part (a).
 To complete the proof it suffices to find a function $f\in
L^p(\R,|x|^{2k}dx)$ so that ${\cal U}(.,t)=P^k_t(f)$ converges
in the $L^p_k$-norm to $f$ and $\|{\cal U}(.,t)\|_{k,p}\leq
B(\alpha,k,t_0,n)A$. Using inequality (\ref{e.V.2"}), we deduce that for
$t\leq t_0$
$$\|{\cal U}(.,t)\|_{k,p}\leq \|{\cal U}(.,t_0)\|_{k,p}+\dint_t^{t_0}\|\partial_s{\cal U}(.,s)\|_{k,p}ds\leq B(n,\alpha,q,t_0)A.$$
On the other hand, by relation (\ref{e.V.2"}), we have $$\|{\cal U}(.,t_1)-{\cal U}(.,t_2)\|_{k,1}\leq
\dint_{t_1}^{t_2}\|\partial_s{\cal U}(.,s)\|_{k,1}ds\leq
B(n,\alpha,q,t_0)A\dint_{t_1}^{t_2}s^{\alpha-1}ds\longrightarrow0,\,\,\mbox{as}\,\,t_1\leq
t_2\longrightarrow0.$$ According to Theorem
\ref{Th.VI.characterization}, there exists $f\in L^p_k$ (it is
uniformly continuous if $p=\infty$) such that ${\cal
U}(x,t)=P^k_tf$. This achieves the proof of the Lemma \ref{Lem.V.2}.
\begin{Rems}\label{Rem.VI.100} We have
\begin{itemize}
  \item By proceeding in same manner as before, we can assert that the Lemma \ref{Lem.V.2} is true when we take ${\cal U}(x,t)$ $k$-temperature on $\R^2_+$
and we replace $k$-Poisson transform by $k$-heat transform.
 \item If $\beta>0$, we define $P^k_t({\cal B}^k_{-\beta})$ as follows
  \begin{equation}\label{e.V.1}
  P_t^k({\cal B}^k_{-\beta})(x)=P^k_t({\cal B}^k_{2-\beta})(x)+\partial^2_tP^k_t({\cal B}^k_{2-\beta})(x),\,\,\,\mbox{when}\,\,\,0<\beta<2,
  \end{equation}
and for arbitrary $\beta>0$ by the rule
$$P_t^k({\cal B}^k_{-\beta})(x)=P^k_{\frac{t}{2}}({\cal B}^k_{-\gamma})\ast_kP^k_{\frac{t}{2}}({\cal B}^k_{-\delta})(x),\,\,\mbox{whenever}\,\,\,\gamma+\delta=\beta.$$
\item If $\beta>0$, we define the $k$-Bessel potential ${\cal J}^k_{-\beta}f(x)$ for a function $f\in L^p(\R,|x|^{2k}dx)$, $1\leq p\leq\infty$, by
$${\cal J}^k_{-\beta}f(x)=\displaystyle\lim_{t\longrightarrow0}P^k_t({\cal B}^k_{-\beta})\ast_kf(x),$$
where the limit is interpreted in $L^p_k$-norm and pointwice a.e.
\end{itemize}
\end{Rems}
\begin{Rem}\label{Rem.V.1}
 For $f\in L^p(\R,|x|^{2k}dx)$, $1\leq p\leq\infty$ and $\beta>0$, the $k$-Poisson transform of ${\cal J}^k_{-\beta}f$,
$P^k_t({\cal J}^k_{-\beta}(f))$, is $k$-harmonic on $\R^2_+$ and
$\|P^k_t({\cal J}^k_{-\beta}(f))\|_{k,p}\leq\|{\cal J}^k_{-\beta}f\|_{k,p}$, for all $t>t_0$, with $t_0>0$.
\end{Rem}
\par We will study the action of the $k$-Bessel potential ${\cal
J}^k_{\beta}$ on the generalized Dunkl-Lipschitz spaces,
$\wedge^k_{\alpha,p,q}(\R)$.
\begin{Th}\label{Th.V.1}
 Let $\alpha>0$, $\beta>0$ and $1\leq p,q\leq\infty$. Then ${\cal
J}^k_{\beta}$ is a topological isomorphism from $\wedge^k_{\alpha,p,q}(\R)$ onto
$\wedge^k_{\alpha+\beta,p,q}(\R)$.
\end{Th}
{\footnotesize\bf{Proof}}\quad If $f\in\wedge^k_{\alpha,p,q}(\R)$, by Lemma
\ref{Lem.V.1}, we have
$$\|{\cal J}^k_{\beta}(f)\|_{\wedge^k_{\alpha+\beta,p,q}}\leq B(k,\beta)\|f\|_{\wedge^k_{\alpha,p,q}}$$
which implies the continuity of ${\cal
J}^k_{\beta}$ from $\wedge^k_{\alpha,p,q}(\R)$ into
$\wedge^k_{\alpha+\beta,p,q}(\R)$. If $f\in\wedge^k_{\alpha+\beta,p,q}(\R)$, we may assume without loss of generality that $\beta\in]0,2[$. Applying the formula (\ref{e.V.1}) and Lemma \ref{Lem.V.1}, we obtain
\begin{equation}\label{e.V.6'}
\|P^k_t({\cal B}^k_{-\beta})\|_{k,1}\leq1+B(k,\beta)t^{-\beta}\leq B(k,\beta),\,\,t\geq1.
\end{equation}
Therefore,
$$\|{\cal J}^k_{-\beta}(P^k_t(f))\|_{k,p}\leq B(k,\beta)\|f\|_{\wedge^k_{\alpha+\beta,p,q}},\,\,t\geq1.$$
From formula (\ref{e.V.6'}) and Proposition \ref{Prop.V.1}(v), a direct verification yields that
$${\cal A}^k_{p,q}(t^{\overline{\alpha}+\overline{\beta}-\alpha}\partial^{\overline{\alpha}+\overline{\beta}}_tP^k_t({\cal J}^k_{-\beta}(f)))\leq
B(k,\alpha,\beta)\|f\|_{\wedge^k_{\alpha+\beta,p,q}}.$$ On the
other hand, by remark \ref{Rem.V.1} and Lemma \ref{Lem.V.2}, there
exists a function $g\in\wedge^k_{\alpha,p,q}(\R)$ satisfying $P^k_t({\cal J}^k_{-\beta}(f))=P^k_t(g)$.
Consequently, we get
$${\cal J}^k_{-\beta}(f)=g\,\,\,\mbox{with}\,\,\,g\in\wedge^k_{\alpha,p,q}(\R)\,\,\,\mbox{and}\,\,\,\|{\cal J}^k_{-\beta}(f)\|
_{\wedge^k_{\alpha,p,q}}\leq
B(k,\alpha,\beta)\|f\|_{\wedge^k_{\alpha+\beta,p,q}}$$ which proves the continuity of ${\cal J}^k_{-\beta}$ from $\wedge^k_{\alpha+\beta,p,q}(\R)$
into $\wedge^k_{\alpha,p,q}(\R)$. We now come to show
 ${\cal J}^k_{-\beta}({\cal
J}^k_{\beta}(f))(x)=f(x)$ a.e., if
$f\in\wedge^k_{\alpha,p,q}(\R)$, $\alpha>0$, which
follows from the fact that\\ $P^k_t({\cal J}^k_{-\beta}({\cal J}^k_{\beta}(f)))(x)=P^k_t(f)(x)$ and similarly, ${\cal J}^k_{\beta}({\cal
J}^k_{-\beta}(f))(x)=f(x)$ a.e., if $f\in\wedge^k_{\alpha+\beta,p,q}(\R)$, $\alpha>0$.
This concludes the proof of the theorem.
 \\\par Before giving a formal definition of the generalized Dunkl-Lipschitz
spaces, we introduce the definition of the space ${\cal L}^p_{\alpha,k}(\R)$.
\begin{Def}
The Lebesgue space
$${\cal L}^p_{\alpha,k}(\R):=\left\{T\in S'(\R):\,\,T={\cal J}^k_{\alpha}(g),\,\,g\in L^p(\R,|x|^{2k}dx)\right\},$$
for $\alpha$ real, $1\leq p\leq\infty$, is called the Dunkl-Sobolev
space of fractional order $\alpha$. Define
$$\|T\|_{k,p,\alpha}:=\|g\|_{k,p}.$$
Thus ${\cal L}^p_{\alpha,k}(\R)$ is a Banach space that is an
isometric image of $L^p(\R,|x|^{2k}dx)$.
\end{Def}
\par Now, following the classical case, see for instance \cite{Taibleson,Flett1}, we are going to define the generalized Dunkl-Lipschitz spaces  $\wedge^k_{\alpha,p,q}(\R)$, for all real $\alpha$.
\begin{Def} Let $p,q\in[1,\infty]$, $\alpha\in\R$ and $n=\overline{(\dfrac{\alpha}{2})}$.
 \\ (i) If $\alpha>0$, $\wedge^k_{\alpha,p,q}(\R)$ is the space of functions of $f\in
 L^p(\R,|x|^{2k}dx)$ for which the $k$-heat transform
 $G^k_t(f)$ of $f$ satisfies the condition that
 $${\cal A}^k_{p,q}(t^{n-\frac{\alpha}{2}}\partial^n_tG^k_t(f))<\infty.$$
 The space is given the norm
 $$\|f\|_{\wedge^k_{\alpha,p,q}}:=\|f\|_{k,p}+{\cal A}^k_{p,q}(t^{n-\frac{\alpha}{2}}\partial^n_tG^k_t(f)).$$
 (ii) If $\alpha\leq0$, $\wedge^k_{\alpha,p,q}(\R)$ is the space of tempered
 distributions $T\in{\cal L}^p_{\alpha-\frac{1}{2},k}(\R)$ for which
 the $k$-heat transform $G^k_t(T)$ of $T$ satisfies the
 condition that
 $${\cal A}^{k,\ast}_{p,q}(t^{n-\frac{\alpha}{2}}\partial^n_tG^k_t(T))<\infty.$$
 The space is given the norm
 $$\|T\|_{\wedge^k_{\alpha,p,q}}:=\|T\|_{k,p,\alpha-\frac{1}{2}}+{\cal A}^{k,\ast}_{p,q}(t^{n-\frac{\alpha}{2}}\partial^n_tG^k_t(T)).$$
\end{Def}
\begin{Lem}\label{Lem.V.2'}
Let $\alpha<0$, $1\leq p\leq\infty$, $T\in {\cal L}^p_{\alpha,k}(\R)$ and let $G^k_t(T)$ be the $k$-heat transform of $T$ on $\R^2_+$.
Then $G^k_t(T)\in{\cal T}^k(\R^2_+)$ and
$$\|G^k_t(T)\|_{k,p}\leq B(k,\alpha)(t^{\frac{1}{2}\alpha}+1)\|T\|_{k,p,\alpha}.$$
\end{Lem}
{\footnotesize\bf{Proof}}\quad From Theorem 3.12 of \cite{N-A-S} and Theorem \ref{Th.IV.3}, the result is proved.
\\\par Now, we want to extend the Theorem \ref{Th.V.1} for all real $\alpha$ and $\beta$. For this, we need the following auxiliary lemmas.
\begin{Lem}\label{Lem.V.3}
Let $H(x,t)$ be absolutely continuous as a
function $t$ for $(x,t)\in\R^2_+$, $t\leq1$. Then for
$\alpha>0$, $p,q\in[1,\infty]$,
$${\cal A}^{k,\ast}_{p,q}(t^{\alpha}H)\leq B(\alpha,q)\left[{\cal A}^{k,\ast}_{p,q}(t^{\alpha+1}\partial_tH)+\|H(.,1)\|_{k,p}\right].$$
\end{Lem}
{\footnotesize\bf{Proof}}\quad We shall prove the Lemma only when $q\in[1,\infty[$, the case $q=\infty$ can be similarly treated. We can write
$$H(x,t)=H(x,1)-\dint_t^1\partial_sH(x,s)ds.$$
From Minkowski's integral inequality, we obtain
$${\cal A}^{k,\ast}_{p,q}(t^{\alpha}H)\leq B(\alpha,q)\|H(.,1)\|_{k,p}+\left\{\dint_0^1\left[t^{\alpha}\dint_t^1
\|\partial_sH(.,s)\|_{k,p}ds\right]^q\dfrac{dt}{t}\right\}^{\frac{1}{q}}.$$
The result announced arises from Hardy inequality.
\begin{Rem}\label{Rem.V.1.}
 Observe that, for $\alpha>0$, the tempered distribution ${\cal B}^k_{\alpha}$ is a
  function in $L^1(\R,|x|^{2k}dx)$. For $\alpha=0$ it is the Dirac
  delta $\delta_0$ and for $-\alpha\in]0,2[$
  $$G^k_t({\cal B}^k_{\alpha})(x)=G^k_t({\cal B}^k_{\alpha+2})(x)-\partial_tG^k_t({\cal B}^k_{\alpha+2})(x)$$
  which is easily verified by taking the Dunkl transform ${\cal F}_k$.
 Similarly, we may construct $G^k_t({\cal B}^k_{\alpha})$ for all
$\alpha<0$ and find in particular that for each $t>0$,
$G^k_t({\cal B}^k_{\alpha})\in L^1(\R,|x|^{2k}dx)$ and is uniformly
bounded in $L^1(\R,|x|^{2k}dx)$ in each proper sub-half space of
$\R^2_+$.
\end{Rem}
\begin{Lem}\label{Lem.V.4}
Let $\alpha$ be real number, $T\in{\cal
L}^p_{\alpha-\frac{1}{2},k}(\R)$ and $n\in\N$,
$n\geq\overline{(\frac{\alpha}{2})}$. Then the norm
$$\|T\|_{k,p,\alpha-\frac{1}{2}}+{\cal A}^{k,\ast}_{p,q}(t^{n-\frac{\alpha}{2}}\partial^n_tG^k_t(T))$$
is equivalent to the norm with $n=\overline{(\frac{\alpha}{2})}$.
\end{Lem}
{\footnotesize\bf{Proof}}\quad If $T\in{\cal
L}^p_{\alpha-\frac{1}{2},k}(\R)$, from Proposition \ref{Prop.II.1}(iv), we have
$$\|\partial^n_tG^k_1(T)\|_{k,p}\leq
B(k,n,\alpha)\|T\|_{k,p,\alpha-\frac{1}{2}},\,\,\,n>l=\overline{(\frac{\alpha}{2})}.$$ Therefore by Lemma \ref{Lem.V.3}, we obtain $${\cal
A}^{k,\ast}_{p,q}(t^{l-\frac{\alpha}{2}}\partial^l_tG^k_t(T))\leq
B(k,\alpha,n)({\cal A}^{k,\ast}_{p,q}(t^{n-\frac{\alpha}{2}}\partial^n_tG^k_t(T))+\|T\|_{k,p,\alpha-\frac{1}{2}}).$$
 Conversely, a direct check shows that
$${\cal A}^{k,\ast}_{p,q}(t^{\beta+1}\partial_tG^k_t(T))\leq B(k,\beta){\cal A}^{k,\ast}_{p,q}(t^{\beta}G^k_t(T)),\,\,\,\beta>0.$$
Thus
$${\cal A}^{k,\ast}_{p,q}(t^{n-\frac{\alpha}{2}}\partial^n_tG^k_t(T))
\leq B(k,\alpha,n){\cal
A}^{k,\ast}_{p,q}(t^{l-\frac{\alpha}{2}}\partial^l_tG^k_t(T)),\,\,\,\mbox{where}\,\,\,n>l=\overline{(\frac{\alpha}{2})},$$
which proves the results.
\begin{Lem}\label{Lem.V.5}
Let $\alpha$ be real, $n=\overline{(\frac{\alpha}{2})}$ and $1\leq p,q\leq\infty$. Then the set of
tempered distributions $T\in{\cal L}^p_{\alpha-\frac{1}{2},k}(\R)$
for which $${\cal
A}^{k,\ast}_{p,q}(t^{n-\frac{\alpha}{2}}\partial^n_tG^k_t(T))<\infty,$$
normed with
\begin{equation}\label{e.V.3}
{\cal
A}^{k,\ast}_{p,q}(t^{n-\frac{\alpha}{2}}\partial^n_tG^k_t(T))+\|T\|_{k,p,\alpha-\frac{1}{2}}
\end{equation} is topologically and algebraically equal to
$\wedge^k_{\alpha,p,q}(\R)$.
\end{Lem}
{\footnotesize\bf{Proof}}\quad By definition of $\wedge^k_{\alpha,p,q}(\R)$, one only needs to consider the case $\alpha>0$.
Assume that $T\in{\cal L}^p_{\alpha-\frac{1}{2},k}(\R)$ and (\ref{e.V.3}) is finite. It is easily seen that
\begin{equation}\label{e.V.4}
{\cal A}^k_{p,q}(t^{n-\frac{\alpha}{2}}\partial^n_tG^k_t(T))\leq
B(k,\alpha,q)\left({\cal
A}^{k,\ast}_{p,q}(t^{n-\frac{\alpha}{2}}\partial^n_tG^k_t(T))+
\|T\|_{k,p,\alpha-\frac{1}{2}}\right),\;\;\alpha>0.
\end{equation}
If $\alpha\geq\frac{1}{2}$, thus
$T\in L^p(\R,|x|^{2k}dx)$ is obvious. On the other hand, if $0<\alpha<\frac{1}{2}$, then for $t\geq1$,
 $\|G^k_t(T)\|_{k,p}\leq B(k,\alpha)\|T\|_{k,p,\alpha-\frac{1}{2}}$. By the relation (\ref{e.V.4}) and Lemma \ref{Lem.V.2},
there exists a function $\psi\in\wedge^k_{\alpha,p,q}(\R)$ such
that $G^k_t(T)=G^k_t(\psi)$ and
$$\|\psi\|_{\wedge^k_{\alpha,p,q}}\leq B(k,\alpha,q)\left\{{\cal A}^{k,\ast}_{p,q}(t^{n-\frac{\alpha}{2}}
\partial^n_tG^k_t(T))+\|T\|_{k,p,\alpha-\frac{1}{2}}\right\}.$$
Now $T$ and $\psi$ have the same $k$-heat transform and thus are
equal as distributions. This implies that $T$ is a function and is in
$L^p(\R,|x|^{2k}dx)$, when $\alpha\in]0,\frac{1}{2}]$. Summarizing, the above two cases show that $T\in\wedge^k_{\alpha,p,q}(\R)$
and
$$\|T\|_{\wedge^k_{\alpha,p,q}}\leq B(k,\alpha,q)\left({\cal A}^{k,\ast}_{p,q}(t^{n-\frac{\alpha}{2}}\partial^n_tG^k_t(T))+
\|T\|_{k,p,\alpha-\frac{1}{2}}\right),\;\;\alpha>0.$$
Conversely, let  $T\in\wedge^k_{\alpha,p,q}(\R)$ and $\|T\|_{\wedge^k_{\alpha,p,q}}$
 is finite. Note that
$${\cal A}^{k,\ast}_{p,q}(t^{n-\frac{\alpha}{2}}\partial^n_tG^k_t(T))\leq{\cal
A}^k_{p,q}(t^{n-\frac{\alpha}{2}}\partial^n_tG^k_t(T))<\infty,\,\,\,\alpha>0.$$
If $\alpha\in]0,\frac{1}{2}]$,then $T\in{\cal L}^p_{\alpha-\frac{1}{2},k}(\R)$ is obvious. If $\alpha>\frac{1}{2}$, thus from Theorem \ref{Th.V.1}, we obtain
$${\cal J}^k_{-(\alpha-\frac{1}{2})}(T)\in\wedge^k_{\frac{1}{2},p,q}(\R)\subset L^p(\R,|x|^{2k}dx)\,\,\,\mbox{and}\,\,\,\|{\cal J}^k_{-(\alpha-\frac{1}{2})}
(T)\|_{k,p}\leq B(k,\alpha)\|T\|_{\wedge^k_{\alpha,p,q}}.$$ Since
$\|T\|_{k,p,\alpha-\frac{1}{2}}=\|{\cal
J}^k_{-(\alpha-\frac{1}{2})}(T)\|_{k,p}$ and
$\|T\|_{\wedge^k_{\alpha,p,q}}$ is finite, the proof is finished.
\begin{Rem}\label{Rem.VI.101}
From Lemmas \ref{Lem.V.1} and \ref{Lem.V.5} for $\beta>0$, Remark
\ref{Rem.V.1.} for $\beta<0$ and Proposition \ref{Prop.II.1}(iv)
for $\beta=0$, we get
$$\|\partial^n_tG^k_t({\cal B}^k_{\beta})\|_{k,1}\leq B(k,\beta) t^{\frac{\beta}{2}-n},\;\;\mbox{where}\;\;n-\frac{\beta}{2}>0\;\;\mbox{and}\;\;t>0.$$
\end{Rem}
\par We can now state the main result of this section.
\begin{Th}\label{Th.V.2}
Let $\alpha$, $\beta$ be real and $1\leq p,q\leq\infty$. Then
${\cal J}^k_{\beta}$ is a topological isomorphism from $\wedge^k_{\alpha,p,q}(\R)$ onto $\wedge^k_{\alpha+\beta,p,q}(\R)$.
\end{Th}
{\footnotesize\bf{Proof}}\quad
Suppose $f\in\wedge^k_{\alpha,p,q}(\R)$, by Remark \ref{Rem.VI.101}, we obtain
$$\|\partial^l_tG^k_t({\cal J}^k_{\beta}(f))\|_{k,p}\leq B(k,\beta)t^{\frac{\beta}{2}-\overline{(\frac{\beta}{2})}}
\|\partial^s_tG^k_{\frac{t}{2}}(f)\|_{k,p},$$ where
$l=\overline{(\frac{\alpha}{2})}+\overline{(\frac{\beta}{2})}$ and
$s=\overline{(\frac{\alpha}{2})}$. As a consequence, we deduce
$${\cal A}^{k,\ast}_{p,q}(t^{l-\frac{\alpha+\beta}{2}}\partial^l_tG^k_t({\cal J}^k_{\beta}(f)))\leq B(k,\beta)
{\cal A}^{k,\ast}_{p,q}(t^{s-\frac{\alpha}{2}}\partial^s_tG^k_t(f)).$$
From Lemmas \ref{Lem.V.4} and \ref{Lem.V.5}, we conclude that
$${\cal J}^k_{\beta}f\in\wedge^k_{\alpha+\beta,p,q}(\R)\,\,\,\mbox{and}\,\,\,\|{\cal J}^k_{\beta}f\|_{\wedge^k_{\alpha+\beta,p,q}}\leq
B(k,\alpha,\beta)\|f\|_{\wedge^k_{\alpha,p,q}}.$$ Moreover, the
following relation
$$G^k_{t_1}({\cal B}^k_{\beta})\ast_kG^k_{t_2}({\cal B}^k_{-\beta})=F^k_{t_1+t_2},\,\,t_1,t_2>0,$$
provide that if $f\in\wedge^k_{\alpha,p,q}(\R)$ then ${\cal
J}^k_{-\beta}({\cal J}^k_{\beta}(f))=f$ as a distribution. Similar conclusions show that if $f\in\wedge^k_{\alpha+\beta,p,q}(\R)$ then
${\cal J}^k_{\beta}({\cal J}^k_{-\beta}(f))=f$ as a distribution. The
announced statement arises.
\begin{Th}\label{Th.V.3}
Let $T\in S'(\R)$. Then for each integer
$n>\overline{(\frac{\alpha}{2})}$ and real number $\beta<\alpha$,
the norm
\begin{equation}\label{e.V.5}
{\cal
A}^{k,\ast}_{p,q}(t^{n-\frac{\alpha}{2}}\partial^n_tG^k_t(T))+\|T\|_{k,p,\beta}
\end{equation}
is equivalent to $\|T\|_{\wedge^k_{\alpha,p,q}}$, where $1\leq
p,q\leq\infty$.
\end{Th}
{\footnotesize\bf{Proof}}\quad Suppose
$T\in\wedge^k_{\alpha,p,q}(\R)$. Since $\alpha-\beta>0$ and by
Theorem \ref{Th.V.2}, we have
$$\|T\|_{k,p,\beta}=\|{\cal J}^k_{-\beta}T\|_{k,p}\leq\|{\cal J}^k_{-\beta}T\|_{\wedge^k_{\alpha-\beta,p,q}}\leq B(k,\alpha,\beta)
\|T\|_{\wedge^k_{\alpha,p,q}}.$$ Then, Lemmas \ref{Lem.V.5} and \ref{Lem.V.4} ensure that relation (\ref{e.V.5}) is finite.\\
Conversely, if relation (\ref{e.V.5}) is finite and let
$l>\overline{(\frac{\alpha-\beta}{2})}$. By Lemmas \ref{Lem.V.3} and \ref{Lem.V.4}, Remark \ref{Rem.VI.101} and change of variables, we have
$${\cal A}^{k,\ast}_{p,q}(t^{l-\frac{\alpha-\beta}{2}}\partial^l_tG^k_t({\cal J}^k_{-\beta}(T)))\leq B(k,n,\alpha,\beta)
\left\{{\cal A}^{k,\ast}_{p,q}(t^{n-\frac{\alpha}{2}}\partial^n_tG^k_t(T))+\|T\|_{k,p,\beta}\right\}$$ and $\|{\cal J}^k_{-\beta}T\|_{k,p}=\|T\|_{k,p,\beta}$.
Note that
$$\|{\cal J}^k_{-\beta}T\|_{\wedge^k_{\alpha-\beta,p,q}}\leq B(k,\alpha,\beta)\left\{{\cal A}^{k,\ast}_{p,q}(t^{l-\frac{\alpha-\beta}{2}}
\partial^l_tG^k_t({\cal J}^k_{-\beta}(T)))+\|{\cal J}^k_{-\beta}T\|_{k,p}\right\},$$
hence from Theorem \ref{Th.V.2}, we obtain
$$\|T\|_{\wedge^k_{\alpha,p,q}}\leq B(k,\alpha,\beta)\|{\cal J}^k_{-\beta}T\|_{\wedge^k_{\alpha-\beta,p,q}}\leq B(k,n,\alpha,\beta)\left\{
{\cal A}^{k,\ast}_{p,q}(t^{n-\frac{\alpha}{2}}\partial^n_tG^k_t(T))+\|T\|_{k,p,\beta}\right\}$$
which prove the theorem.\\\\
{\footnotesize{\bf{Note}}}
\par We are essentially defining $\wedge^k_{-\alpha,p,q}(\R)$ to
be ${\cal
J}^k_{-\alpha-\frac{1}{2}}(\wedge^k_{\frac{1}{2},p,q}(\R))$, $\alpha>0$. The
choice of $\frac{1}{2}$ is arbitrary. Any $\beta>0$, would work as
well.
\\\par The remainder of this section is devoted to some properties and embedding theorems for the spaces $\wedge^k_{\alpha,p,q}(\R)$.
\begin{Th}\label{Th.V.4}
Let $f$ in
$\wedge^k_{\alpha_0,p_0,q_0}(\R)\cap\wedge^k_{\alpha_1,p_1,q_1}(\R)$,
then $f$ belongs to $\wedge^k_{\alpha,p,q}(\R)$ and we have
$$\|f\|_{\wedge^k_{\alpha,p,q}}\leq
B(k,\alpha_0,\alpha_1)\|f\|^{1-\theta}_{\wedge^k_{\alpha_0,p_0,q_0}}\|f\|^{\theta}_{\wedge_{\alpha_1,p_1,q_1}},$$
where $\alpha=(1-\theta)\alpha_0+\theta\alpha_1$,
$\dfrac{1}{p}=\dfrac{1-\theta}{p_0}+\dfrac{\theta}{p_1}$,
$\dfrac{1}{q}=\dfrac{1-\theta}{q_0}+\dfrac{\theta}{q_1}$, and
$\theta\in[0,1]$. In particular
\\ (a)
$\|f\|_{k,p,\beta}\leq\|f\|^{1-\theta}_{k,p_0,\beta}\|f\|^{\theta}_{k,p_1,\beta}$,
$\beta<\min(\alpha_0,\alpha_1)$.
\\ (b) ${\cal A}^k_{p,q}(t^{n-\frac{\alpha}{2}}\partial^n_tG^k_t(f))\leq\left[{\cal A}^k_{p_0,q_0}(t^{n-\frac{\alpha_0}{2}}\partial^n_tG^k_t(f))
\right]^{1-\theta}\left[{\cal
A}^k_{p_1,q_1}(t^{n-\frac{\alpha_1}{2}}\partial^n_tG^k_t(f))
\right]^{\theta}$, where
$n>\max(\frac{\alpha_0}{2},\frac{\alpha_1}{2})$.
\end{Th}
{\footnotesize\bf{Proof}}\quad This can be proved from Theorem
\ref{Th.V.3} and the Logarithmic convexity of the $L^p_k$-norms.
\\\par Let us study some inclusions among the generalized Dunkl-Lipschitz spaces :
\begin{Lem}\label{Lem.V.6}
The continuous embedding
  $$\wedge^k_{\alpha_1,p,q_1}(\R)\hookrightarrow\wedge^k_{\alpha_2,p,q_2}(\R)$$
  holds if either
  \\ (i) if $\alpha_1>\alpha_2$ ( then $q_1$ and $q_2$ need not be
  related), or
  \\ (ii) if $\alpha_1=\alpha_2$ and $q_1\leq q_2$.
\end{Lem}
{\footnotesize\bf{Proof}}\quad We give the argument for
$q\neq\infty$. The case $q=\infty$ is done similarly. We may
suppose $0<\alpha_2<\alpha_1<1$. Let
$f\in\wedge^k_{\alpha_1,p,q_1}(\R)$ and consider first the case
$q_1=q_2$. In the one hand, it is easily to see that
$${\cal
A}^{k,\ast}_{p,q_1}(t^{1-\alpha_2}\partial_tP^k_t(f))\leq\|f\|_{\wedge^k_{\alpha_1,p,q_1}}.$$
In the other hand, using the fact that $\|\partial_tP^k_t(f)\|_{k,p}\leq B(k)t^{-1}\|f\|_{k,p}$, we get
$$\left\{\dint_1^{\infty}\left[t^{1-\alpha_2}\|\partial_tP^k_t(f)\|_{k,p}\right]^{q_1}\dfrac{dt}{t}\right\}^{\frac{1}{q_1}}
\leq B(k,\alpha_2,q_1)\|f\|_{\wedge^k_{\alpha_1,p,q_1}}$$ which
proves that $\wedge^k_{\alpha_1,p,q_2}(\R)\hookrightarrow\wedge^k_{\alpha_2,p,q_2}(\R)$. Moreover, if $q_1<q_2$, Lemma 5.2 of \cite{S-K} and
Lemma 1.2 of \cite{Johnson} show that
$\wedge^k_{\alpha_1,p,q_1}(\R)\hookrightarrow\wedge^k_{\alpha_1,p,q_2}(\R)$.
Hence
$\wedge^k_{\alpha_1,p,q_1}(\R)\hookrightarrow\wedge^k_{\alpha_1,p,q_2}(\R)\hookrightarrow\wedge^k_{\alpha_2,p,q_2}(\R)$.
If $q_1>q_2$, let $\frac{1}{s}=\frac{1}{q_2}-\frac{1}{q_1}$. Applying
H\"older's inequality and analogous reasoning as before finish the
proof of the lemma.
\begin{Lem}\label{Lem.V.7}
If $1\leq p_1\leq p_2$ and
$\alpha_1-\frac{2k+1}{p_1}=\alpha_2-\frac{2k+1}{p_2}$, we have the
continuous embedding
$$\wedge^k_{\alpha_1,p_1,q}(\R)\hookrightarrow\wedge^k_{\alpha_2,p_2,q}(\R).$$
\end{Lem}
{\footnotesize\bf{Proof}}\quad We may assume that
$0<\alpha_1,\alpha_2<1$. If $f\in\wedge^k_{\alpha_1,p_1,q}(\R)$, Young's inequality yields that
$$\|\partial_tP^k_t(f)\|_{k,p_2}\leq\|\partial_tP^k_{\frac{t}{2}}(f)\|_{k,p_1}\|P^k_{\frac{t}{2}}\|_{k,s}\leq
B(k,p_1,p_2)
t^{(-\frac{1}{p_1}+\frac{1}{p_2})(2k+1)}\|\partial_tP^k_{\frac{t}{2}}(f)\|_{k,p_1},
$$ where $\frac{1}{s}=\frac{1}{p_2}-\frac{1}{p_1}+1$. Hence ${\cal A}^k_{p_2,q}(t^{1-\alpha_2}\partial_tP^k_t(f))\leq B(k,\alpha_1,p_1,p_2){\cal
A}^k_{p_1,q}(t^{1-\alpha_1}\partial_tP^k_t(f))$. On the other
hand, for $t\geq1$, $\|P^k_t(f)\|_{k,p_2}\leq
B(k,p_1,p_2)\|f\|_{k,p_1}$ and therefore by Lemma \ref{Lem.V.2},
we can deduce that $f\in\wedge^k_{\alpha_2,p_2,q}(\R)$ and
$\|f\|_{\wedge^k_{\alpha_2,p_2,q}}\leq
B(k,\alpha_1,p_1,p_2)\|f\|_{\wedge^k_{\alpha_1,p_1,q}}$ which end the proof.
\\\par As consequence of Lemmas \ref{Lem.V.6} and
\ref{Lem.V.7}, we deduce the following theorem :
\begin{Th}\label{Th.V.5}
Let $\alpha_1,\alpha_2\in\R$ and $1\leq p_1\leq p_2\leq\infty$, then we have the continuous
embedding
$$\wedge^k_{\alpha_1,p_1,q_1}(\R)\hookrightarrow\wedge^k_{\alpha_2,p_2,q_2}(\R)$$
if $\alpha_1-\frac{2k+1}{p_1}>\alpha_2-\frac{2k+1}{p_2}$ or if
$\alpha_1-\frac{2k+1}{p_1}=\alpha_2-\frac{2k+1}{p_2}$ and $1\leq
q_1\leq q_2\leq\infty$.
\end{Th}
\par The action of Dunkl derivatives on Dunkl-Lipschitz spaces is
as follows :
\begin{Prop}\label{Prop.V.2}
Let $\alpha>0$, $1\leq p,q\leq\infty$ and $0\leq n\leq\alpha$.
Then the norm $\|f\|_{k,p}+\|{\cal
D}^n_kf\|_{\wedge^k_{\alpha-n,p,q}}$ is equivalent to
$\|f\|_{\wedge^k_{\alpha,p,q}}$.
\end{Prop}
{\footnotesize\bf{Proof}}\quad If $\|f\|_{\wedge^k_{k,p,q}}$ is
finite, then according to the Proposition \ref{Prop.V.1}(v) and Remark (5.14) of \cite{S-K}, it is easy to see that
$${\cal A}^k_{p,q}(t^{\overline{\alpha}-(\alpha-n)}\partial_t^{\overline{\alpha}}{\cal D}^n_kP^k_t(f))\leq B(k,\alpha,n)\|f\|
_{\wedge^k_{\alpha,p,q}},$$ and $$\|{\cal
D}^n_kP^k_t(f)\|_{k,p}\leq B(k,n)\|f\|_{k,p},\,\,\,t\geq1.$$ Thus by
Lemma \ref{Lem.V.2}, we deduce
that there exists $g\in\wedge^k_{\alpha-n,p,q}(\R)$ such that
${\cal D}^n_kP^k_t(f)=P^k_t(g)$ and
$\|g\|_{\wedge^k_{\alpha-n,p,q}}\leq
B(k,\alpha,n)\|f\|_{\wedge^k_{\alpha,p,q}}$. On the other hand,
since ${\cal D}^n_kP^k_t(f)=P^k_t({\cal D}^n_kf)$ (in the distribution sense), we have $P^k_t(g)=P^k_t({\cal D}^n_kf)$. Letting $t\longrightarrow0$ yields that $g={\cal D}^n_kf$. An easy check shows the converse result.
\begin{Lem}\label{Lem.V.8}
If $f\in\wedge^k_{\alpha,\infty,q}(\R)$, $\alpha\in]0,1[$, then
$f$ is uniformly continuous.
\end{Lem}
{\footnotesize\bf{Proof}}\quad It suffices to show that
$\|\triangle_{y,k}f\|_{k,\infty}\rightarrow0$ as $y\rightarrow0$.
By Theorem \ref{Th.V.5},
$f\in\wedge^k_{\alpha,\infty,\infty}(\R)$, so
$\|\triangle_{y,k}f\|_{k,\infty}\leq A|y|^{\alpha}$ and thus tends
to zero as $y\rightarrow0$.
\begin{Th}\label{Th.V.6}
$\wedge^k_{\alpha,p,q}(\R)$ is complete if $1\leq p,q\leq\infty$ and $\alpha\in\R$.
\end{Th}
{\footnotesize\bf{Proof}}\quad By Theorem \ref{Th.V.2}, we may suppose $\alpha\in]0,1[$. If $(f_n)$ is a Cauchy sequence in
$\wedge^k_{\alpha,p,q}(\R)$, then $(f_n)$ is obviously Cauchy sequence in
$L^p_k$, and therefore converges in $L^p_k$ to a function $f$. Hence $\|\partial_tP^k_t(f_s)\|_{k,p}\rightarrow\|\partial_tP^k_t(f)\|_{k,p}$ as
$s\rightarrow\infty$ and for $m=1,2,\cdots$, $\|\partial_t(P^k_tf_m-P^k_tf_s)\|_{k,p}\rightarrow\|\partial_t(P^k_tf_m-P^k_tf)\|_{k,p}$ as $s\rightarrow\infty$.
Consequently, by Fatou's Lemma, we have
$${\cal A}^k_{p,q}(t^{1-\alpha}\partial_t(P^k_tf_m-P^k_tf))\leq\epsilon_m=\displaystyle\lim_{s\rightarrow\infty}\inf{\cal A}^k_{p,q}(t^{1-\alpha}
\partial_t(P^k_tf_m-P^k_tf_s))_{\overrightarrow{m\rightarrow\infty}}0,$$ and ${\cal A}^k_{p,q}(t^{1-\alpha}\partial_tP^k_t(f))
\leq\displaystyle\lim_{s\rightarrow\infty}\inf\|f_s\|_{\wedge^k_{\alpha,p,q}}<\infty$.
So that $f\in\wedge^k_{\alpha,p,q}(\R)$ and $f_m \rightarrow f$, as $m\rightarrow\infty$, in
$\wedge^k_{\alpha,p,q}(\R)$ which conclude the proof.
\\\par The object of the next section will be to derive a similar result for $k$-temperatures on $\R^2_+$.
\section{\footnotesize{{\bf{Dunkl-Lipschitz Spaces of $k$-Temperatures}}}}
We shall define a generalized Dunkl-Lipschitz space of
$k$-temperatures on $\R^2_+$ which will be denoted by ${\cal T}\wedge^k_{\alpha,p,q}(\R^2_+)$ and prove that various norms are equivalent to
our original definition. Finally, the isomorphism of ${\cal T}\wedge^k_{\alpha,p,q}(\R^2_+)$ and $\wedge^k_{\alpha,p,q}(\R)$ is established.
\\ We begin this section by stating the following standard Lemmas.
\begin{Def}
Let $\alpha$ be a real number. For any $k$-temperature ${\cal U}$ in
${\cal T}^k(\R^2_+)$, $1\leq p\leq\infty$ and $1\leq
q\leq\infty$, let
$${\cal E}^{k,\alpha}_{p,q}({\cal U}):=\left\{
\begin{array}{rcl}
\left\{\dint_0^{+\infty}t^{q-1}e^{-t}\|{\cal
J}^k_{-\alpha-2}{\cal U}(.,t)\|^q_{k,p}dt\right\}^{\frac{1}{q}}&(1\leq
q<\infty),\\
\displaystyle\sup_{t>0}\left\{te^{-t}\|{\cal
J}^k_{-\alpha-2}{\cal U}(.,t)\|_{k,p}\right\}\,\,\,\,\,\,\,&(q=\infty),
\end{array}
\right.
$$
with infinite values being allowed.
\end{Def}
\begin{Lem}\label{Lem.VI.1}
Let $\alpha$, ${\cal U}$, $p$, $q$ be as in the above definition and let
$\gamma$ be a real number. Then $${\cal
E}^{k,\alpha}_{p,q}({\cal U})={\cal E}^{k,\alpha+\gamma}_{p,q}({\cal J
}^k_{\gamma}{\cal U}).$$
\end{Lem}
{\footnotesize\bf{Proof}}\quad By Theorem \ref{Th.IV.2}, ${\cal
J}^k_{-\alpha-2}{\cal U}={\cal J}^k_{-\alpha-\gamma-2}({\cal
J}^k_{\gamma}{\cal U})$ which implies that ${\cal
E}^{k,\alpha}_{p,q}({\cal U})={\cal E}^{k,\alpha+\gamma}_{p,q}({\cal J
}^k_{\gamma}{\cal U})$.
\begin{Def}\label{Def.VII.102}
Let $1\leq p,q\leq\infty$, let $\alpha$, $\beta$ be real numbers such
that $\beta>\alpha$. For any $k$-temperature ${\cal U}$ in ${\cal
T}^k(\R^2_+)$, let
$${\cal E}^{k,\alpha,\beta}_{p,q}({\cal U}):=\left\{
\begin{array}{rcl}
\left\{\dint_0^{+\infty}t^{\frac{1}{2}q(\beta-\alpha)-1}e^{-t}\|{\cal
J}^k_{-\beta}{\cal U}(.,t)\|^q_{k,p}dt\right\}^{\frac{1}{q}}&(1\leq
q<\infty),\\
\displaystyle\sup_{t>0}\left\{t^{\frac{1}{2}(\beta-\alpha)}e^{-t}\|{\cal
J}^k_{-\beta}{\cal U}(.,t)\|_{k,p}\right\}\,\,\,\,\,\,\,&(q=\infty),
\end{array}
\right.
$$
and $${\cal L}^k_p({\cal U}):=\displaystyle\sup_{t\geq\frac{1}{2}}\|{\cal U}(.,t)\|_{k,p}.$$
\end{Def}
\begin{Rem}\label{Rem.VII.1}
Let $1\leq p,q\leq\infty$, and $\gamma$ be real number. If ${\cal U}\in{\cal T}^k(\R^2_+)$ and ${\cal E}^{k,\alpha}_{p,q}({\cal U})<\infty$, where $\alpha$ is real, so that Theorem \ref{Th.IV.4} and Corollary \ref{Cor.IV.1} yield that for each $a>0$ there exists a positive constant $B$ such that for all $t\geq a$
$$\|{\cal J}^k_{\gamma}{\cal U}(.,t)\|_{k,p}\leq B(k,\alpha,\gamma,q,a){\cal E}^{k,\alpha}_{p,q}({\cal U}).$$
\end{Rem}
\begin{Lem}\label{Lem.VI.2}
Let $\alpha$, $\beta$, ${\cal U}$, $p$, $q$ be as in
definition \ref{Def.VII.102}. Then
\\ (i) ${\cal E}^{k,\alpha}_{p,q}({\cal U})$ is equivalent to
${\cal E}^{k,\alpha,\beta}_{p,q}({\cal U})$.
\\ (ii) ${\cal E}^{k,\alpha}_{p,q}({\cal U})$ is equivalent to ${\cal A}^{k,\ast}_{p,q}\left(t^{\frac{1}{2}(\beta-\alpha)}{\cal J}^k_{-\beta}{\cal U}\right)+{\cal L}^k_p({\cal U})$.
\end{Lem}
{\footnotesize\bf{Proof}}\quad The proof is a simple consequence of Remark \ref{Rem.VII.1}, Theorem \ref{Th.IV.5} and Corollary \ref{Cor.IV.1}.
\begin{Lem}\label{Lem.VI.3}
Let $\alpha$ be real number, ${\cal U}\in{\cal T}^k(\R^2_+)$, $1\leq
p\leq\infty$, $1\leq q\leq\infty$, and $n$ be a non-negative
integer greater than $\frac{\alpha}{2}$. Then ${\cal
E}^{k,\alpha}_{p,q}({\cal U})$ is equivalent to ${\cal A}^{k,\ast}_{p,q}\left(t^{n-\frac{1}{2}\alpha}\partial^n_t{\cal U}\right)+{\cal L}^k_p({\cal U})$.
\end{Lem}
{\footnotesize\bf{Proof}}\quad If $n=0$, the result will be obtained from Lemma \ref{Lem.VI.2}(ii).
First suppose that ${\cal E}^{k,\alpha}_{p,q}({\cal U})<\infty$. For $i=0,1,\cdots,n-1$, we have
$$\|{\cal J}^k_{-2i}{\cal U}(.,t)\|_{k,p}\leq\|{\cal J}^k_{-2n}{\cal U}(.,t)\|_{k,p}$$ and since $\partial^n_t{\cal U}(.,t)$ is a linear combination of ${\cal U}(.,t),\;\;{\cal J}^k_{-2}{\cal U}(.,t),\;\;\cdots,{\cal J}^k_{-2n}{\cal U}(.,t)$, it follows that
\begin{equation}\label{e.VI.1}
\|\partial^n_t{\cal U}(.,t)\|_{k,p}\leq B(k,n)\|{\cal J}^k_{-2n}{\cal U}(.,t)\|_{k,p}
\end{equation}
and therefore by Lemma \ref{Lem.VI.2}(ii), we obtained
$${\cal L}^k_p({\cal U})+{\cal A}^{k,\ast}_{p,q}\left(t^{n-\frac{1}{2}\alpha}\partial^n_t{\cal U}\right)\leq{\cal L}^k_p({\cal U})+{\cal A}^{k,\ast}_{p,q}\left(t^{n-\frac{1}{2}\alpha}{\cal J}^k_{-2n}{\cal U}\right)\leq B(k,n,\alpha,q){\cal E}^{k,\alpha}_{p,q}(U).
$$ Conversely, suppose ${\cal L}^k_p({\cal U})+{\cal A}^{k,\ast}_{p,q}\left(t^{n-\frac{1}{2}\alpha}\partial^n_t{\cal U}\right)$. From Theorem \ref{Th.III.2}, Minkowski's integral inequality, relation (\ref{e.I.4}) and Proposition \ref{Prop.II.1}(iv),
we deduce that for $i=1,2\cdots n$
$$\displaystyle\sup_{t\geq1}\|\partial^i_t{\cal U}(.,t)\|_{k,p}\leq B(k,i){\cal L}^k_p({\cal U})$$ and
\begin{equation}\label{e.VI.2}
\|\partial^i_t{\cal U}(.,t)\|_{k,p}\leq B(k,n) {\cal L}^k_p({\cal U})+\|\partial^n_t{\cal U}(.,t)\|_{k,p}.
\end{equation}
Thus $${\cal A}^{k,\ast}_{p,q}\left(t^{n-\frac{1}{2}\alpha}{\cal J}^k_{-2n}{\cal U}\right)\leq B(k,n,\alpha,q)\left({\cal L}^k_p({\cal U})+{\cal A}^{k,\ast}_{p,q}\left(t^{n-\frac{1}{2}\alpha}\partial^n_t{\cal U}\right)\right).$$
Again Lemma \ref{Lem.VI.2}(ii) shows the desired result.
\\\par Now we turn to the definitions of the generalized Dunkl-Lipschitz space of $k$-temperatures on $\R^2_+$.
\begin{Def}
Let $\alpha$ be a real number, $1\leq p\leq\infty$, $1\leq
q\leq\infty$. We define
$${\cal T}\wedge^k_{\alpha,p,q}(\R^2_+):=\left\{{\cal U}\in{\cal T}^k(\R^2_+)\,:{\cal E}^{k,\alpha}_{p,q}({\cal U})<\infty\right\};$$
$${\cal T}\lambda^k_{\alpha,p,\infty}(\R^2_+):=\left\{{\cal U}\in{\cal T}\wedge^k_{\alpha,p,\infty}(\R^2_+)\,:\|{\cal J}^k_{-\alpha-2}{\cal U}(.,t)\|_{k,p}=
\circ(t^{-1})\,\,\,\mbox{as}\,\,\,t\longrightarrow0^+\right\}.$$
Then, ${\cal E}^{k,\alpha}_{p,q}$ is a norm on ${\cal T}\wedge^k_{\alpha,p,q}(\R^2_+)$.
\end{Def}
\par First we give :
\begin{Lem}\label{Lem.VI'.2}
Let $1\leq p,q\leq\infty$, $\alpha$ and $\gamma$ be real
numbers. Then ${\cal J}^k_{\gamma}$ is an isometric isomorphism of
${\cal T}\wedge^k_{\alpha,p,q}(\R^2_+)$ (${\cal
T}\lambda^k_{\alpha,p,\infty}(\R^2_+)$ resp.) onto ${\cal
T}\wedge^k_{\alpha+\gamma,p,q}(\R^2_+)$ (${\cal
T}\lambda^k_{\alpha+\gamma,p,\infty}(\R^2_+)$ resp.) with inverse
${\cal J}^k_{-\gamma}$.
\end{Lem}
{\footnotesize\bf{Proof}}\quad Since ${\cal J}^k_{-\alpha-2}{\cal U}={\cal J}^k_{-\alpha-\gamma-2}\left({\cal J}^k_{\gamma}{\cal U}\right)$, then Corollary \ref{Cor.IV.01}
proves the result.
\\\\ The basic properties of the spaces ${\cal
T}\wedge^k_{\alpha,p,q}(\R^2_+)$ lie in the following theorem :
\begin{Th}\label{Th.VI'.1}
Let $1\leq p,q\leq\infty$ and $\alpha$ be a real number.
\\ (i) If $1\leq q_1< q_2<\infty$, we have the continuous embedding
$${\cal T}\wedge^k_{\alpha,p,q_1}(\R^2_+)\hookrightarrow{\cal T}\wedge^k_{\alpha,p,q_2}(\R^2_+)\hookrightarrow{\cal
T}\lambda^k_{\alpha,p,\infty}(\R^2_+)\hookrightarrow{\cal
T}\wedge^k_{\alpha,p,\infty}(\R^2_+).$$
\\ (ii) If $\beta$ is a real number such that
$\beta>\alpha$, then ${\cal E}^{k,\alpha,\beta}_{p,q}$ is an
equivalent norm on ${\cal T}\wedge^k_{\alpha,p,q}(\R^2_+)$; moreover
${\cal U}\in{\cal T}\lambda^k_{\alpha,p,\infty}(\R^2_+)$ if and only if
${\cal U}\in{\cal T}\wedge^k_{\alpha,p,\infty}(\R^2_+)$ and $\|{\cal
J}^k_{-\beta}{\cal U}(.,t)\|_{k,p}=\circ(t^{-\frac{1}{2}(\beta-\alpha)})$
as $t\longrightarrow0^+$.
\\ (iii) If $n$ is a non-negative integer greater than
$\frac{1}{2}\alpha$, then ${\cal A}^{k,\ast}_{p,q}\left(t^{n-\frac{1}{2}\alpha}\partial^n_t{\cal U}\right)+{\cal L}^k_p({\cal U})$
is an equivalent norm on ${\cal T}\wedge^k_{\alpha,p,q}(\R^2_+)$.
\\ (iv) The spaces ${\cal T}\wedge^k_{\alpha,p,q}(\R^2_+)$,
where $p$, $q$ are fixed and $\alpha$ varies, are isomorphic to
one another. The same conclusion holds for the spaces ${\cal
T}\lambda^k_{\alpha,p,\infty}(\R^2_+)$.
\end{Th}
{\footnotesize\bf{Proof}}\quad (i) follows easily from Theorem \ref{Th.IV.4}.
(ii) is an easy consequence of Lemma \ref{Lem.VI.2} and Theorem
\ref{Th.IV.5}(iii). (iii) is derived from Lemma \ref{Lem.VI.3}. To prove (iv), let $\delta$ be another real number. It
then follows from Lemma \ref{Lem.VI'.2} that ${\cal J}^k_{-n}$ is
an isometric isomorphism of ${\cal T}\wedge^k_{\delta,p,q}(\R^2_+)$ (${\cal T}\lambda^k_{\delta,p,\infty}(\R^2_+)$ resp.) onto ${\cal
T}\wedge^k_{\delta-n,p,q}(\R^2_+)$ (${\cal
T}\lambda^k_{\delta-n,p,\infty}(\R^2_+)$ resp.); denote its inverse by
$({\cal J}^k_{-n})^{-1}$. This Lemma again implies that ${\cal
J}^k_{\delta-\alpha-n}$ is an isometric isomorphism of ${\cal
T}\wedge^k_{\alpha,p,q}(\R^2_+)$ (${\cal
T}\lambda^k_{\alpha,p,\infty}(\R^2_+)$ resp.) onto ${\cal
T}\wedge^k_{\delta-n,p,q}(\R^2_+)$ (${\cal
T}\lambda^k_{\delta-n,p,\infty}(\R^2_+)$ resp.). Consequently, $({\cal
J}^k_{-n})^{-1}\circ{\cal J}^k_{\delta-\alpha-n}$ is an isometric
isomorphism of ${\cal T}\wedge^k_{\alpha,p,q}(\R^2_+)$ (${\cal
T}\lambda^k_{\alpha,p,\infty}(\R^2_+)$ resp.) onto ${\cal
T}\wedge^k_{\delta,p,q}(\R^2_+)$ (${\cal
T}\lambda^k_{\delta,p,\infty}(\R^2_+)$ resp.).
\\\par The following theorem establish the relation between $\wedge^k_{\alpha,p,q}(\R)$ and ${\cal T}\wedge^k_{\alpha,p,q}(\R^2_+)$.
\begin{Th}\label{Th.VII.1} If $1\leq p,q\leq\infty$ and $\alpha$ is real, then the $k$-heat transform is a topological isomorphism from $\wedge^k_{\alpha,p,q}(\R)$ onto ${\cal T}\wedge^k_{\alpha,p,q}(\R^2_+)$. Moreover if $f\in\wedge^k_{\alpha,p,q}(\R)$, then $G^k_t(f)\in{\cal T}\wedge^k_{\alpha,p,q}(\R^2_+)$ and ${\cal E}^{k,\alpha}_{p,q}(G^k_t(f))\leq B(k,\alpha)\|f\|_{\wedge^k_{\alpha,p,q}}$. Conversely, if ${\cal U}\in{\cal T}\wedge^k_{\alpha,p,q}(\R^2_+)$, then there exists $f\in\wedge^k_{\alpha,p,q}(\R)$ such that $${\cal U}(.,t)=G^k_t(f)(.),\;\;t>0,\;\;\mbox{and}\;\;\|f\|_{\wedge^k_{\alpha,p,q}}\leq B(k,\alpha){\cal E}^{k,\alpha}_{p,q}({\cal U}).$$
\end{Th}
{\footnotesize\bf{Proof}}\quad Let $f\in\wedge^k_{\alpha,p,q}(\R)$, by Theorem \ref{Th.II.2}, Lemmas \ref{Lem.V.2'} and \ref{Lem.VI.3}, we deduce that $$G^k_t(f)\in{\cal T}\wedge^k_{\alpha,p,q}(\R)\;\;\mbox{and}\;\;{\cal E}^{k,\alpha}_{p,q}(G^k_t(f))\leq B(k,\alpha)\|f\|_{\wedge^k_{\alpha,p,q}}.$$ To prove the converse  we proceed first in case $\alpha>0$. For ${\cal U}\in{\cal T}\wedge^k_{\alpha,p,q}(\R^2_+)$, let ${\cal V}(.,t)={\cal J}^k_{-\alpha-2}{\cal U}(.,t)$, $t>0$, then for $s>0$
$${\cal U}(x,s)=\dfrac{1}{\Gamma(\frac{\alpha}{2}+1)}\dint_0^{+\infty}\xi^{\frac{\alpha}{2}}e^{-\xi}{\cal V}(x,\xi+s)d\xi.$$
Moreover, by Theorem \ref{Th.IV.4} yields
\begin{equation}\label{e.VII.1}
\|{\cal J}^k_{-\alpha-2}{\cal U}(.,t)\|_{k,p}\leq B(q)(t^{-1}+1){\cal E}^{k,\alpha}_{p,q}({\cal U})
\end{equation}
which together with Minkowski's integral inequality, we find that
$$\|{\cal U}(.,s)\|_{k,p}\leq B(q,\alpha){\cal E}^{k,\alpha}_{p,q}({\cal U})\dint_0^{+\infty}\xi^{\frac{\alpha}{2}}e^{-\xi}(\xi^{-1}+1)d\xi=B(q,\alpha){\cal E}^{k,\alpha}_{p,q}({\cal U}),\,\,\,\mbox{if}\,\,\,1\leq p\leq\infty.$$
On the one hand, for $p=1$ and $\epsilon>0$, from inequality (\ref{e.VII.1}), we can find $\delta$ satisfying $0<\delta<1$ such that $\|{\cal V}(.,t)\|_{k,1}\leq\epsilon t^{-1-\frac{1}{4}\alpha}$ for $0< t\leq\delta$. On the other hand, by a simple verification yields $\|{\cal U}(.,s)-{\cal U}(.,s')\|_{k,1}\rightarrow0$ as $s,s'\rightarrow0$. Summarizing the above two cases show that from Remark \ref{Rem.VI.Temperature}, there exists a function $f\in L^p(\R,|x|^{2k}dx)$, $1\leq p\leq\infty$, such that ${\cal U}(.,t)=G^k_t(f)(.)$. Next, in case $\alpha\leq0$, then using Lemma \ref{Lem.VI'.2}, we have
$${\cal J}^k_{-\alpha+\frac{1}{2}}{\cal U}\in{\cal T}\wedge^k_{\frac{1}{2},p,q}(\R^2_+)\;\;\mbox{and}\;\;{\cal E}^{k,\frac{1}{2}}_{p,q}({\cal J}^k_{-\alpha+\frac{1}{2}}{\cal U})\leq B{\cal E}^{k,\alpha}_{p,q}({\cal U}).$$ Applying the above case $\alpha>0$, then there exists $g\in L^p(\R,|x|^{2k}dx)$, $p\in[1,\infty]$, such that ${\cal J}^k_{-\alpha+\frac{1}{2}}{\cal U}(.,t)=G^k_t(g)(.)$, and $\|g\|_{k,p}\leq B{\cal E}^{k,\alpha}_{p,q}({\cal U})$. Due to Theorem 3.12 for \cite{N-A-S},
$${\cal U}(.,t)=G^k_t(f)(.),\;\;f={\cal J}^k_{\alpha-\frac{1}{2}}(g)\;\;\mbox{and}\;\;\|f\|_{k,p,\alpha-\frac{1}{2}}=\|g\|_{k,p}\leq B{\cal E}^{k,\alpha}_{p,q}({\cal U}).$$
By Proposition \ref{Prop.II.1}(iv), we obtain for $\alpha>0$ $${\cal A}^k_{p,q}(t^{n-\frac{\alpha}{2}}\partial^n_t{\cal U})\leq {\cal A}^{k,\ast}_{p,q}\left(t^{n-\frac{1}{2}\alpha}\partial^n_t{\cal U}\right)+B(k,\alpha){\cal L}^k_p({\cal U}),\,\,n=\overline{(\frac{\alpha}{2})}.$$ Therefore by Lemma \ref{Lem.VI.3}, we obtain $\|f\|_{\wedge^k_{\alpha,p,q}}\leq B{\cal E}^{k,\alpha}_{p,q}({\cal U})$, $\alpha\in\R$, and the theorem is proved.
\begin{Th}\label{Th.VI'.2}
Let $1\leq p<r\leq\infty$, $1\leq q\leq\infty$, $\alpha$ be a real
number and $\delta=\frac{1}{p}-\frac{1}{r}$. Then
$$(i)\,\,\,{\cal T}\wedge^k_{\alpha,p,q}(\R^2_+)\hookrightarrow{\cal
T}\wedge^k_{\alpha-\delta(2k+1),r,q}(\R^2_+), \,\,\,(ii)\,\,\,{\cal
T}\lambda^k_{\alpha,p,\infty}(\R^2_+)\hookrightarrow{\cal
T}\lambda^k_{\alpha-\delta(2k+1),r,\infty}(\R^2_+).$$
\end{Th}
{\footnotesize\bf{Proof}}\quad Let $h$ such that $\frac{1}{r}=\frac{1}{p}+\frac{1}{h}-1$, ($\frac{1}{h}=1-\delta$). We give the argument for $q\neq\infty$. The case $q=\infty$ is done similarly. Let ${\cal U}$ be in ${\cal T}\wedge^k_{\alpha,p,q}(\R^2_+)$ and $\beta$ be a real number greater than $\alpha$. Theorem \ref{Th.VI'.1}(ii) implies that ${\cal E}^{k,\alpha,\beta}_{p,q}({\cal U})$ is equivalent to ${\cal E}^{k,\alpha}_{p,q}({\cal U})$, $\beta>\alpha$. Then $t\mapsto\|{\cal J}^k_{-\beta}{\cal U}(.,t)\|_{k,p}$ is locally integrable on $]0,\infty[$, so the semi-group formula holds for ${\cal J}^k_{-\beta}{\cal U}$. By Theorem \ref{Th.III.2} and Young's inequality (Proposition 7.2 of \cite{Xu1}), we have
$$\|{\cal J}^k_{-\beta}{\cal U}(.,t)\|_{k,r}\leq\|{\cal J}^k_{-\beta}{\cal U}(.,\frac{t}{2})\|_{k,p}\|F^k_{\frac{t}{2}}\|_{k,h}.$$
By a simple verification, we deduce that $\|F^k_{\frac{t}{2}}\|_{k,h}\leq B(k,p,r)t^{-(k+\frac{1}{2})\delta}$.
Hence, we obtain
$$\|{\cal J}^k_{-\beta}{\cal U}(.,t)\|_{k,r}\leq B(k,p,r)t^{-(k+\frac{1}{2})\delta}\|{\cal J}^k_{-\beta}{\cal U}(.,\frac{t}{2})\|_{k,p}.$$
Therefore
$${\cal E}^{k,\alpha-2k\delta-\delta,\beta}_{r,q}({\cal U})=\left\{\dint_0^{+\infty}t^{\frac{1}{2}q(\beta-\alpha+2k\delta+\delta)-1}e^{-t}\|{\cal J}^k_{-\beta}{\cal U}(.,t)\|^q_{k,r}dt\right\}^{\frac{1}{q}}$$$$\;\;\;\leq B(k,p,r)\left\{\dint_0^{+\infty}t^{\frac{1}{2}q(\beta-\alpha)-1}e^{-t}\|{\cal J}^k_{-\beta}{\cal U}(.,\frac{t}{2})\|^q_{k,p}dt\right\}^{\frac{1}{q}}\leq B(k,p,\alpha,\beta,r){\cal E}^{k,\alpha,\beta}_{p,q}({\cal U}),$$
from which we obtain the part (i) after making use of Theorem \ref{Th.VI'.1}(ii) again. We proceed in the same way to prove the assertion (ii).
\begin{Rem}
In view of the isometry between $\wedge^k_{\alpha,p,q}(\R)$ and ${\cal T}\wedge^k_{\alpha,p,q}(\R^2_+)$, the same result of Theorem \ref{Th.VI'.2} holds for spaces $\wedge^k_{\alpha,p,q}(\R)$.
\end{Rem}


{\footnotesize
\begin{thebibliography}{99}
\bibitem{Ch-An-Sa-Si}Abdelkefi, C., Anker, J.P., Sassi, F., Sifi, M.:
Besov-type spaces on $\R^d$ and integrability for the Dunkl
transform. Symmetry Integrability Geom. Methods Appl. {\bf{5}}, Paper 019, 15 pp. (2009)
\bibitem{Ch-Sa}Abdelkefi, C., Sassi, F.: Various characterizations of Besov-Dunkl spaces. Int. J.
Pure Appl. Math. {\bf{39}}(4), 475-488 (2007)
\bibitem{Ch-Si}Abdelkefi, C., Sifi, M.: Characterization of Besov spaces for the Dunkl operator on
the real line. J. Inequal. Pure Appl. Math. {\bf{8}}(3), Article 73, 11 pp. (2007)
\bibitem{Aron}Aronszajn, N., Smith, K.T.: Theory of Bessel potentials.
I. Ann. de  L'Inst. Fourier {\bf{11}}, 385-475(1961)
\bibitem{Bateman}Bateman, H.: Tables of integral transforms. Vol II, New York  (1954)
\bibitem{N-A-S}Ben Salem, N., El Garna, A., Kallel, S.: Bessel and Flett potentials associated with Dunkl operators
on $\R^d$. Methods and Application of Analysis {\bf{15}}(4), 477-494 (2008)
\bibitem{N-S}Ben Salem, N., Kallel, S.:
Mean-periodic functions associated with the Dunkl
operators. Integral Transforms Spec. Funct. {\bf{15}}(2), 155-179 (2004)
\bibitem{BLA}Bouguila, R., Lazhari, M.N., Assal, M.: Besov spaces associated with Dunkl's operator. Integral
Transforms Spec. Funct. {\bf{18}}(8), 545-557 (2007)
\bibitem{Jeu1}de Jeu, M.F.E.: The Dunkl transform. Invent. Math. {\bf{113}}, 147-162 (1993)
\bibitem{Dunkl3}Dunkl, C.F.: Differential-difference operators
associated to reflection groups. Trans. Am. Math. Soc. {\bf{311}}, 167-183 (1989)
\bibitem{Dunkl1}Dunkl, C.F.: Integral kernels with reflection
group invariance. Can. J. Math. {\bf{43}}, 1213-1227 (1991)
\bibitem{Dunkl2}Dunkl, C.F.: Hankel transforms associated to finite reflection
groups. Contemp. Math. {\bf{138}}, 123-138 (1992)
\bibitem{Flett1}Flett, T.M.: Temperatures, Bessel potentials and Lipschitz
 spaces. Proc. London Math. Soc. {\bf{22}}(3), 385-451 (1971)
\bibitem{Friedman}Friedman, A.: Partial differential equations of parabolic type (Englewood Cliffs) (1964)
\bibitem{Zh-Ji}Li, Zh., Liao, J. : Harmonic analysis associated with the one-dimensional Dunkl transform. Constr. Approx. {\bf{37}}, 233-281 (2013)
\bibitem{Johnson}Johnson, R.: Temperatures, Riesz potentials, and the Lipschitz spaces of Hertz.
 Proc. London Math. Soc. {\bf{27}}(3), 290-316 (1973)
\bibitem{S-K}Kallel, S.: Characterization of function spaces for the Dunkl type operator on the real line. Potential anal. DOI 10.1007/s11118-013-9366-5 (2013)
\bibitem{Lot1}Kamoun, L.: Besov-type spaces for the Dunkl operator on the real line. J. Comput. Appl. Math. {\bf{199}}, 56-67 (2007)
\bibitem{Lot2}Kamoun, L., Nagzaoui, S.: Lipschitz spaces associated with reflection group $\Z^d_2$. Commun.
Math. Anal. {\bf{7}}(1), 21-36 (2009)
\bibitem{Mourou}Mourou, M.A.: Taylor series associated with a differential-difference
operator on the real line. J. Comput. Appl. Math. {\bf{153}}, 343-354 (2003)
\bibitem{Ros}R\"osler, M.: Bessel-type signed hypergroups on
$\R$. In : Heyer, H., Mukherjea, A. (eds.) Probability Measures on Groups and Related Structures
XI (Oberwolfach, 1994), pp. 292-304. World Sci.
Publ. (1995)
\bibitem{Ros-Voit}R\"osler, M., Voit, M.: Dunkl theory, convolution algebras, and related Markov processes.
 In: Harmonic and stochastic analysis of Dunkl processes; eds. P. Graczyk et al. Travaux en cours {\bf{71}},
 Hermann, Paris, pp. 1-112. Preprint version (2008)
\bibitem{Taibleson}Taibleson, M.H.: On the theory of Lipschitz spaces of distributions on Euclidean
$n$-spaces. I. Principal properties, J. Math. Mech. {\bf{13}}(3), 407-479 (1964)
\bibitem{Xu1}Thangavelu, S., Xu, Y.: Convolution operator and maximal function
for the Dunkl transform. J. Anal. Math. {\bf{97}}, 25-55 (2005)
\bibitem{Xu2}Thangavelu, S., Xu, Y.: Riesz transforms and Riesz potentials
for the Dunkl transform. J. Comput. Appl. Math. {\bf{199}}, 181-195 (2007)
\bibitem{Trim2}Trim\`eche, K.: The Dunkl intertwining operator on spaces of functions and
 distributions and integral representations of its dual.
 Integral Transforms Spec. Funct. {\bf{12}}, 349-374 (2001)

\end{thebibliography}

}
 \end{document}